\documentclass[11pt]{article}
\usepackage{geometry}                
\usepackage{amsmath,amssymb}
\usepackage{amsthm}
\usepackage{hyperref}
\usepackage{graphicx,epsfig,color}
\usepackage{booktabs,bm,multirow}
\usepackage{cases}
\usepackage[caption=false]{subfig}
\usepackage{lscape} 
\usepackage{float}
\newtheorem{theorem}{Theorem}

\newtheorem{corollary}[theorem]{Corollary}

\newtheorem{lemma}[theorem]{Lemma}

\def\mbbr{\mathbb R}
\def\mcalk{\mathcal K}

\def\l{\left} 
\def\r{\right} 
\def\mbf{\mathbf}
\def\rma{{\rm A}}
\def\rmd{{\rm D}}
\def\rmr{{\rm R}}
\def\spa{{\rm span}}
\def\rank{{\rm rank}}
\def\ran{{\rm range}}
\def\nul{{\rm null}} 
\def\ind{{\rm ind}} 

\DeclareMathOperator*{\argmin}{argmin}

\def\dsp{\displaystyle} 
\def\nn{\nonumber}

\def\wt{\widetilde}
\def\wh{\widehat}

\newcommand{\beq}{\begin{equation}}\newcommand{\eeq}{\end{equation}}
\newcommand{\bit}{\begin{itemize}}\newcommand{\eit}{\end{itemize}}
\newcommand{\bem}{\begin{bmatrix}}\newcommand{\eem}{\end{bmatrix}}

\title{Obtaining the pseudoinverse solution of singular range-symmetric linear systems with GMRES-type methods\thanks{This work was supported by the National Natural Science Foundation of China (No.12171403 and No.11771364), the Natural Science Foundation of Fujian Province of China (No.2020J01030), and the Fundamental Research Funds for the Central Universities (No.20720210032).}}

\author{Kui Du\thanks{Corresponding author. School of Mathematical Sciences and Fujian Provincial Key Laboratory of Mathematical Modeling and High Performance Scientific Computing, Xiamen University, Xiamen 361005, China (kuidu@xmu.edu.cn).},\quad Jia-Jun Fan\thanks{School of Mathematical Sciences, Xiamen University, Xiamen 361005, China ({fanjiajun@stu.xmu.edu.cn}).},\quad Fang Wang\thanks{School of Mathematical Sciences, Xiamen University, Xiamen 361005, China ({fangwang@stu.xmu.edu.cn}).}} 
\date{}                                           

\begin{document}
\maketitle

\begin{abstract} It is well known that for singular inconsistent range-symmetric linear systems, the generalized minimal residual (GMRES) method determines a least squares solution without breakdown. The reached least squares solution may be or not be the pseudoinverse solution. We show that a lift strategy can be used to obtain the pseudoinverse solution. In addition, we propose a new iterative method named RSMAR (minimum $\mbf A$-residual) for range-symmetric linear systems $\mbf A\mbf x=\mbf b$. At step $k$ RSMAR minimizes $\|\mbf A\mbf r_k\|$ in the $k$th Krylov subspace generated with $\{\mbf A, \mbf r_0\}$ rather than $\|\mbf r_k\|$, where $\mbf r_k$ is the $k$th residual vector and $\|\cdot\|$ denotes the Euclidean vector norm. We show that RSMAR and GMRES terminate with the same least squares solution when applied to range-symmetric linear systems. 
We provide two implementations for RSMAR. Our numerical experiments show that RSMAR is the most suitable method among GMRES-type methods for singular inconsistent range-symmetric linear systems. \

\vspace{2mm} 

\vspace{2mm}  
\end{abstract}

\section{Introduction} 

We consider the linear system of equations \beq\label{lin}\mbf A\mbf x=\mbf b,\eeq where $\mbf b\in\mbbr^n$ is a  vector, and $\mbf A\in\mbbr^{n\times n}$ is a large singular range-symmetric (i.e., $\ran(\mbf A)=\ran(\mbf A^\top)$) matrix for which matrix-vector products $\mbf A\mbf v$ can be computed efficiently for any vector $\mbf v\in\mbbr^n$. For any $\mbf b\in\mbbr^n$, we seek the unique solution $\mbf x_\star$ that solves the problem \beq\label{normal} \min{\|\mbf x\|}\qquad\mbox{subject to}\qquad \mbf A^\top\mbf A\mbf x=\mbf A^\top\mbf b.\eeq It is clear that $\mbf x_\star$ is the unique minimum Euclidean norm solution to \eqref{lin} if $\mbf b\in\ran(\mbf A)$ and the unique minimum Euclidean norm least squares solution otherwise. Here we call $\mbf x_\star$ the pseudoinvese solution of \eqref{lin}. 

If system \eqref{lin} is consistent (i.e., $\mbf b\in\ran(\mbf A)$), then the GMRES method by Saad and Schulz \cite{saad1986gmres} determines the pseudoinvese solution without breakdown. If system \eqref{lin} is inconsistent (i.e., $\mbf b\notin\ran(\mbf A)$), then GMRES determines a least squares solution without breakdown, and the reached least squares solution may be or not be the pseudoinverse solution. We refer to \cite[section 2]{brown1997gmres} for the above statements. Applicable solvers for the pseudoinverse solution of \eqref{lin} with arbitrary $\mbf b\in\mbbr^n$ would be the RRGMRES method by Calvetti, Lewis, and Reichel \cite{calvetti2000gmres} and the DGMRES method by Sidi \cite{sidi2001dgmre}. Let $\mbf x_k$ be an approximate solution to $\mbf x_\star$ with residual $\mbf r_k=\mbf b-\mbf A\mbf x_k$. Assume that $\mbf x_0=\mbf 0$. At step $k$, GMRES minimizes $\|\mbf r_k\|$ over the $k$th Krylov subspace $$\mcalk_k(\mbf A,\mbf b):=\spa\{\mbf b,\mbf A\mbf b,\ldots,\mbf A^{k-1}\mbf b\},$$ RRGMRES minimizes $\|\mbf r_k\|$ over the $k$th Krylov subspace $\mcalk_k(\mbf A,\mbf A\mbf b)$ (which belongs to $\ran(\mbf A)$ and thus is called range-restricted), and DGMRES minimizes $\|\mbf A^\alpha\mbf r_k\|$ over $\mcalk_k(\mbf A,\mbf A\mbf b)$. Here, $\alpha$ is the index of $\mbf A$, the size of a largest Jordan block associated with zero eigenvalue. If $\mbf A$ is range-symmetric, then $\alpha=1$ (see section \ref{pre}). 

When $\mbf A$ is symmetric, GMRES is theoretically equivalent to MINRES \cite{paige1975solut}. Hence, MINRES determines the pseudoinverse solution if $\mbf b\in\ran(\mbf A)$ and a least squares solution (but not necessarily the pseudoinverse solution) otherwise. MINRES-QLP \cite{choi2011minre}, a variant of MINRES, is an applicable solver for the pseudoinverse solution of \eqref{lin} with symmetric $\mbf A$. On ill-conditioned symmetric linear systems (singular or not), MINRES-QLP can give more accurate solutions than MINRES. We mention that Liu, Milzarek, and Roosta \cite{liu2023obtai} proposed a novel and remarkably simple lifting strategy for MINRES to obtain the pseudoinverse solution when $\mbf b\notin\ran(\mbf A)$. The lifting strategy seamlessly integrates with the final MINRES iteration. Compared to MINRES-QLP, the lifted MINRES method can obtain the pseudoinverse solution with negligible additional computational costs. 

Recently, Montoison, Orban, and Saunders \cite{montoison2023minar} proposed an iterative method, named MINARES, for solving symmetric linear systems. At step $k$, MINARES minimizes $\|\mbf A\mbf r_k\|$ over the $k$th Krylov subspace $\mcalk_k(\mbf A,\mbf b)$. Their numerical experiments with MINRES-QLP \cite{choi2011minre} and LSMR \cite{fong2011lsmr} show that MINARES is the most suitable Krylov method for inconsistent symmetric linear systems. Like MINRES, MINARES determines the pseudoinverse solution if $\mbf b\in\ran(\mbf A)$ and a least squares solution (but not necessarily the pseudoinverse solution) otherwise. 

In this paper, we consider GMRES-type methods for range-symmetric linear systems. We mainly focus on the singular case and seek the pesudoinverse solution. 

The main contributions of this work are as follows. (i) We show that the lifting strategy in \cite{liu2023obtai} also works for GMRES on singular inconsistent range-symmetric linear systems (see Theorem \ref{main1}). (ii) We propose a new Krylov subspace method called RSMAR (Range-Symmetric Minimum $\mbf A$-Residual) for computing a solution to range-symmetric linear systems. At step $k$, RSMAR minimizes $\|\mbf A\mbf r_k\|$ over the $k$th Krylov subspace $\mcalk_k(\mbf A,\mbf b)$, and thus is theoretically equivalent to MINARES when applied to symmetric linear systems. (iii) We show that RSMAR and GMRES terminate with the same least squares solution for range-symmetric linear systems, which implies that MINARES and MINRES also terminates with the same least squares solution for symmetric linear systems. (iv) We propose two implementations for RSMAR, named RSMAR-I and RSMAR-II. RSMAR-I is inspired by the implementation for the simpler GMRES method \cite{walker1994simpl}, and RSMAR-II is inspired by the implementation of RRGMRES \cite{neuman2012algor,neuman2012imple}. The MINARES implementation in \cite[section 4]{montoison2023minar}  can be viewed as a short recurrence variant of RSMAR-II. We provide a new implementation for MINARES, which can be viewed as a short recurrence variant of RSMAR-I. (v) Our numerical experiments show that RSMAR-II is the preferable algorithm for singular inconsistent range-symmetric linear systems. 

The paper is organized as follows. In the rest of this section, we give other related research. In section 2, we provide clarification of notation, some properties of the Moore--Penrose inverse and the Drazin inverse, and some useful results for Krylov subspaces. In section 3, we consider four GMRES-type methods (GMRES, RRGMRES, RSMAR, and DGMRES) for singular range-symmetric linear systems, prove our main theoretical results, and provide two implementations for RSMAR. In section 4, we consider two MINRES-type methods (MINRES and MINARES), and provide a new implementation for MINARES. In section 5, some numerical experiments are performed to compare the performance of the methods considered in this paper. Finally, we give some concluding remarks and possible future work in section 6.

{\it Other related research}. In addition to \cite{brown1997gmres}, there exist numerous studies on GMRES for singular linear systems in the literature; see, for example, \cite{ipsen1998idea,smoch1999some,calvetti2000gmres,cao2002note,reichel2005break,smoch2007spect,du2008inexa,zhang2010note,hayami2011geome,morikuni2018gmres,sugihara2023gmres}. GMRES on almost singular (or numerically singular) systems was analyzed in \cite{elden2012solvi}. GMRES for least squares problems was discussed in \cite{hayami2010gmres,morikuni2013inner,morikuni2015conve}. Some convergence properties of Krylov subspace methods for singular linear systems with arbitrary index were discussed in \cite{wei2000conve}. The $m$-shift GMRES method (for which RRGMRES is a special case) was proposed and studied in \cite{bellalij2015some}. Stagnation analysis, restart variant, and convergence rate of DGMRES were studied in \cite{zhou2004stagn}, \cite{zhou2006analy}, and \cite{greenbaum2018conve}, respectively. A simpler DGMRES was proposed in \cite{zhou2010simpl}. For singular symmetric linear systems, some preconditioning techniques for MINRES were considered in \cite{sugihara2020right,hong2022preco}. 

\section{Preliminaries}\label{pre}

\subsection{Notation} Lowercase (uppercase) boldface letters are reserved for column vectors (matrices). Lowercase lightface letters are reserved for scalars. For any vector $\mbf v\in\mbbr^n$, we use $\mbf v^\top$ and $\|\mbf v\|$ to denote the transpose and the Euclidean norm of $\mbf v$, respectively. We use $\mbf I_k$ to denote the $k\times k$ identity matrix, and use $\mbf e_i$ to denote the $i$th column of the identity matrix $\mbf I$ whose order is clear from the context. We use $\mbf 0$ to denote the zero vector (or matrix) of appropriate size. For any matrix $\mbf M\in\mbbr^{n\times n}$, we use $\mbf M^\top$, $\mbf M^\dag$, $\mbf M^\rmd$, and $\|\mbf M\|$ to denote the transpose, the Moore--Penrose inverse, the Drazin inverse, and the spectral norm of $\mbf M$, respectively. For nonsingular $\mbf M$, we use $\mbf M^{-1}$ to denote its inverse. We denote the null space and range of $\mbf M$ by $\nul(\mbf M)$ and $\ran(\mbf M)$, respectively. For a matrix $\mbf M$, its condition number is denoted by $\kappa(\mbf M)=\|\mbf M\|\|\mbf M^\dag\|$, which is the ratio of the largest singular value of $\mbf M$ to the smallest positive one. Throughout the paper, we assume that exact arithmetic is used for all theoretical discussions. 

\subsection{Pseudoinverse solution and Drazin-inverse solution}

The Moore--Penrose inverse of $\mbf A$ is defined as the unique matrix $\mbf A^\dag$ satisfying $$\mbf A\mbf A^\dag\mbf A=\mbf A,\quad \mbf A^\dag\mbf A\mbf A^\dag=\mbf A^\dag,\quad (\mbf A\mbf A^\dag)^\top=\mbf A\mbf A^\dag,\quad (\mbf A^\dag\mbf A)^\top=\mbf A^\dag\mbf A.$$ 
If $\mbf A$ has a zero eigenvalue with index $\alpha$ (the size of a largest Jordan block associated with zero eigenvalue, also called the index of $\mbf A$, denoted by $\ind(\mbf A)$), then the Drazin inverse of $\mbf A$ is defined as the unique matrix $\mbf A^\rmd$ satisfying $$\mbf A^\rmd\mbf A\mbf A^\rmd=\mbf A^\rmd,\quad \mbf A^\rmd\mbf A=\mbf A\mbf A^\rmd, \quad \mbf A^{\alpha+1}\mbf A^\rmd=\mbf A^\alpha.$$  The Drazin inverse $\mbf A^\rmd$ is always expressible as a polynomial of $\mbf A$. We refer to \cite{ben2003gener,wang2018gener} for more properties of the Moore--Penrose inverse and the Drazin inverse. The vector $\mbf A^\dag\mbf b$ is called the pseudoinverse solution of $\mbf A\mbf x=\mbf b$, and the vector $\mbf x^\rmd=\mbf A^\rmd\mbf b$ is called the Drazin inverse solution. The unique solution $\mbf x_\star$ of \eqref{normal} is $\mbf A^\dag\mbf b$.

Let $\mbf x_0\in\mbbr^n$ be a given vector. It is clear that  the vector $\mbf A^\dag\mbf b+(\mbf I-\mbf A^\dag\mbf A)\mbf x_0$ is the orthogonal projection of $\mbf x_0$ onto the solution set $\{\mbf x\in\mbbr^n\ |\ \mbf A\mbf x=\mbf b\}$ if $\mbf b\in\ran(\mbf A)$, and onto the least squares solution set $\{\mbf x\in\mbbr^n\ |\ \mbf A^\top\mbf A\mbf x=\mbf A^\top\mbf b\}$ if $\mbf b\notin\ran(\mbf A)$.


\subsection{Range-symmetric matrix}

A matrix $\mbf A$ is called range-symmetric if $\ran(\mbf A)=\ran(\mbf A^\top)$. 
A range-symmetric matrix $\mbf A$ can be expressed as (see, for example, \cite[Theorem 2.5]{hayami2011geome}) $$\mbf A=\mbf U\bem \mbf C & \mbf 0 \\ \mbf 0 & \mbf 0\\  \eem\mbf U^\top,$$  where the matrix $\mbf C$ is invertible, and the matrix $\mbf U$ is orthogonal. In this case, we have $$\mbf A^\dag=\mbf A^\rmd=\mbf U\bem \mbf C^{-1} & \mbf 0 \\ \mbf 0 & \mbf 0\\  \eem\mbf U^\top.$$ It is clear that range-symmetric $\mbf A$ has index one. When $\mbf A$ is range-symmetric, the linear system $\mbf A^2\mbf x=\mbf A\mbf b$ and the normal equations $\mbf A^\top\mbf A\mbf x=\mbf A^\top\mbf b$ have the same solution set, i.e., the affine set $\mbf A^\dag\mbf b+\nul(\mbf A)$.

\subsection{Krylov subspaces}
 Beginning with an initial approximate solution $\mbf x_0$, at step $k$ a Krylov subspace method \cite{ipsen1998idea} for solving \eqref{lin} generates an approximate solution $\mbf x_k\in\mbf x_0+\mcalk_k(\mbf A,\mbf r_0)$, where $\mbf r_0:=\mbf b-\mbf A\mbf x_0$ and $\mcalk_k(\mbf A,\mbf r_0)$ is the $k$th Krylov subspace $$\mcalk_k(\mbf A,\mbf r_0):=\spa\{\mbf r_0,\mbf A\mbf r_0,\ldots,\mbf A^{k-1}\mbf r_0\}.$$ 
It is well known (see, for example, \cite{reichel2005break}) that there exists an integer $\ell$ satisfying $$\dim\mcalk_k(\mbf A,\mbf r_0)=\begin{cases} k & \mbox{if }  k\leq \ell \\ \ell & \mbox{if } k\geq\ell+1.\end{cases}$$ We know that $\ell$ is the maximal dimension of Krylov subspace generated with the matrix-vector pair $\{\mbf A,\mbf r_0\}$.
The Arnoldi process \cite{arnoldi1951princ} with the matrix-vector pair $\{\mbf A,\mbf r_0\}$ constructs a sequence of orthonormal vectors $\{\mbf v_k\}$ such that $\mbf v_1=\mbf r_0/\beta_1$ with $\beta_1=\|\mbf r_0\|$, $\mbf V_k^\top\mbf V_k=\mbf I_k$, and $$\mbf A\mbf V_k=\mbf V_{k+1}\mbf H_{k+1,k},$$
where $\mbf V_k:=\bem\mbf v_1 & \mbf v_2 & \cdots & \mbf v_k\eem$, and $$\mbf H_{k+1,k}:=\bem h_{11} & \cdots & h_{1k}\\ h_{21} &\ddots & \vdots \\ & \ddots & h_{kk}\\ && h_{k+1,k}\eem$$ is a $(k+1)\times k$ upper-Hessenberg matrix. Let $\mbf H_k$ denote the leading $k\times k$ submatrix of $\mbf H_{k+1,k}$. We have $\mbf H_k=\mbf V_k^\top\mbf A\mbf V_k$. The Arnoldi process with $\{\mbf A,\mbf r_0\}$ terminates at step $\ell$ with $h_{\ell+1,\ell}=0$ and $h_{k+1,k}>0$ for each $1\leq k\leq\ell-1$. We have $\rank(\mbf H_{\ell,\ell-1})=\ell-1$ and \beq\label{ell}\mbf A\mbf V_\ell=\mbf V_\ell\mbf H_\ell.\eeq 
The first $k\leq\ell$ columns of $\mbf V_\ell$ form an orthonormal basis of $\mcalk_k(\mbf A,\mbf r_0)$. We have the following estimates on the number of least squares solution in the affine space $\mbf x_0+\mcalk_{\ell-1}(\mbf A,\mbf r_0)$, and on the number of solution in the affine space $\mbf x_0+\mcalk_\ell(\mbf A,\mbf r_0)$.

\begin{theorem}\label{one} There is at most one least squares solution in $\mbf x_0+\mcalk_{\ell-1}(\mbf A,\mbf r_0)$ if $\mbf b\notin\ran(\mbf A)$, and at most one solution in $\mbf x_0+\mcalk_\ell(\mbf A,\mbf r_0)$ if $\mbf b\in\ran(\mbf A)$.
\end{theorem}
\proof 
Assume that $\mbf x\in\mbf x_0+\mcalk_{\ell-1}(\mbf A,\mbf r_0)$ and $\mbf y\in\mbf x_0+\mcalk_{\ell-1}(\mbf A,\mbf r_0)$ are two least squares solutions. Then we have $\mbf x-\mbf y\in\nul(\mbf A)\cap\mcalk_{\ell-1}(\mbf A,\mbf r_0)$. This means there exists a vector $\mbf z\in\mbbr^{\ell-1}$ such that $$\mbf x-\mbf y =\mbf V_{\ell-1}\mbf z,\quad \mbf A\mbf V_{\ell-1}\mbf z=\mbf V_\ell\mbf H_{\ell,\ell-1}\mbf z=\mbf 0.$$ Thus, $\mbf z=\mbf 0$, which implies $\mbf x=\mbf y$.

The second part is a direct result of Ipsen and Meyer \cite{ipsen1998idea}. If $\mbf b\in\ran(\mbf A^\alpha)$, then the unique solution is $\mbf x_0+\mbf A^\rmd\mbf r_0=\mbf A^\rmd\mbf b+(\mbf I-\mbf A^\rmd\mbf A)\mbf x_0\in\mbf x_0+\mcalk_\ell(\mbf A,\mbf r_0)$. If $\mbf b\notin\ran(\mbf A^\alpha)$, then no solution lies in $\mbf x_0+\mcalk_\ell(\mbf A,\mbf r_0)$. 
\endproof

Now we give some existing results about the matrix $\mbf H_\ell$ in \eqref{ell}. If $\mbf H_\ell$ is nonsingular, then by $$\mbf b-\mbf A(\mbf x_0+\beta_1\mbf V_\ell\mbf H_\ell^{-1}\mbf e_1)=\mbf r_0-\beta_1\mbf A\mbf V_\ell\mbf H_\ell^{-1}\mbf e_1=\mbf r_0-\beta_1\mbf V_\ell\mbf e_1=\mbf 0,$$ we have $\mbf b\in\ran(\mbf A)$. Hence, if $\mbf b\notin\ran(\mbf A)$, then $\mbf H_\ell$ must be singular and $\rank(\mbf H_\ell)=\ell-1$ (because $\mbf H_\ell$ has a nonsingular $(\ell-1)\times (\ell-1)$ upper triangular submatrix). If $\ran(\mbf A)=\ran(\mbf A^\top)$ and $\mbf b\in\ran(\mbf A)$, then $\mbf H_\ell$ must be nonsingular (see, for example, \cite{brown1997gmres}).

The Krylov subspace $\mcalk_k(\mbf A,\mbf A\mbf r_0)=\mbf A\mcalk_k(\mbf A,\mbf r_0)$ is important in our analysis. Using $\dim\mcalk_k(\mbf A,\mbf r_0)=k$ for $k\leq\ell$, $\dim\mcalk_k(\mbf A,\mbf r_0)=\ell$ for $k>\ell$, and $$\mcalk_k(\mbf A,\mbf A\mbf r_0)=\spa\{\mbf A\mbf r_0,\mbf A^2\mbf r_0,\ldots,\mbf A^k\mbf r_0\}\subseteq\mcalk_{k+1}(\mbf A,\mbf r_0),$$ 
we have $\dim\mcalk_k(\mbf A,\mbf A\mbf r_0)=k$ for each $1\leq k\leq \ell-1$ and $\dim\mcalk_k(\mbf A,\mbf A\mbf r_0)\leq\ell$ for all $k\geq\ell$. Using \eqref{ell}, we further have $\dim\mcalk_\ell(\mbf A,\mbf A\mbf r_0)=\ell$ if $\mbf H_\ell$ is nonsingular and $\dim\mcalk_\ell(\mbf A,\mbf A\mbf r_0)=\dim\mcalk_{\ell-1}(\mbf A,\mbf A\mbf r_0)=\ell-1$ otherwise. Let $m$ denote the maximal dimension of Krylov subspace generated with $\{\mbf A,\mbf A\mbf r_0\}$. We have \beq\label{mm}m=\dim\mcalk_\ell(\mbf A,\mbf A\mbf r_0)=\begin{cases} \ell & \mbox{if }  \mbf H_\ell \mbox{ is nonsingular} \\ \ell-1 & \mbox{if } \mbf H_\ell \mbox{ is singular}.\end{cases}\eeq  
Using the Arnoldi process with the matrix-vector pair $\{\mbf A,\mbf A\mbf r_0\}$, we obtain an orthonormal basis, denoted by $\{\wh{\mbf v}_k\}$, for $\mcalk_k(\mbf A,\mbf A\mbf r_0)$ such that $\wh{\mbf v}_1=\mbf A\mbf r_0/\wh\beta_1$ with $\wh\beta_1=\|\mbf A\mbf r_0\|$, $\wh{\mbf V}_k^\top\wh{\mbf V}_k=\mbf I_k$, and $$\mbf A\wh{\mbf V}_k=\wh{\mbf V}_{k+1}\wh{\mbf H}_{k+1,k},$$
where $\wh{\mbf V}_k:=\bem\wh{\mbf v}_1 & \wh{\mbf v}_2 & \cdots & \wh{\mbf v}_k\eem$, and $$\wh{\mbf H}_{k+1,k}:=\bem \wh h_{11} & \cdots & \wh h_{1k}\\ \wh h_{21} &\ddots & \vdots \\ & \ddots & \wh h_{kk}\\ && \wh h_{k+1,k}\eem.$$ 
Let $\wh{\mbf H}_k$ denote the leading $k\times k$ submatrix of $\wh{\mbf H}_{k+1,k}$. We have $\wh{\mbf H}_k=\wh{\mbf V}_k^\top\mbf A\wh{\mbf V}_k$. The Arnoldi process with $\{\mbf A,\mbf A\mbf r_0\}$ terminates at step $m$ with $\wh h_{m+1,m}=0$ and $\wh h_{k+1,k}>0$ for each $1\leq k\leq m-1$. We have $\rank(\wh{\mbf H}_{m,m-1})=m-1$ and $$\mbf A\wh{\mbf V}_m=\wh{\mbf V}_m\wh{\mbf H}_m.$$ The first $k\leq m$ columns of $\wh{\mbf V}_m$ form an orthonormal basis of $\mcalk_k(\mbf A,\mbf A\mbf r_0)$. In the following theorem, we give a condition ensuring the invertibility of the matrix $\wh{\mbf H}_m$.

\begin{theorem}\label{index1} If the index of $\mbf A$ is one $(\ind(\mbf A)=1)$, then the matrix $\wh{\mbf H}_m$ is nonsingular.
\end{theorem}
\proof Ipsen and Mayer \cite[Theorem 2]{ipsen1998idea} proved that $\mbf A\mbf x=\mbf b$ has a Krylov solution in $\mcalk_n(\mbf A,\mbf b)$ if and only if $\mbf b\in\ran(\mbf A^\alpha)$, where $\alpha=\ind(\mbf A)$. Thus, if $\ind(\mbf A)=1$, then by $\mbf A\mbf r_0\in\ran(\mbf A)$ we know that $\mbf A\mbf y=\mbf A\mbf r_0$ has a Krylov solution in $\mcalk_m(\mbf A,\mbf A\mbf r_0)$. That is to say there exists a vector $\mbf z\in\mbbr^m$ such that $\mbf A\wh{\mbf V}_m\mbf z=\mbf A\mbf r_0$. Using $\mbf A\wh{\mbf V}_m=\wh{\mbf V}_m\wh{\mbf H}_m$ and $\mbf A\mbf r_0=\wh\beta_1\wh{\mbf v}_1$, we have $ \wh{\mbf V}_m\wh{\mbf H}_m\mbf z=\wh\beta_1\wh{\mbf V}_m\mbf e_1$. This means $\wh{\mbf H}_m\mbf z=\wh\beta_1\mbf e_1$ is consistent. Therefore, we have $\rank(\wh{\mbf H}_m)= \rank(\bem \wh\beta_1\mbf e_1 & \wh{\mbf H}_m\eem)=m$ (since the matrix consisting of the first $m$ columns of $\bem \wh\beta_1\mbf e_1 & \wh{\mbf H}_m\eem$ is nonsingular upper triangular). Hence, $\wh{\mbf H}_m$ is nonsingular.    
\endproof

Note that range-symmetric $\mbf A$ has index one. A direct result of Theorem \ref{index1} is that $\wh{\mbf H}_m$ is nonsingular if $\mbf A$ is range-symmetric.  

\subsection{Summary of some scalars, vectors, and matrices}

For clarity in the following discussions, we list the frequently used scalars, vectors, and matrices in this paper in the following table.

\begin{table}[H]
\caption{Frequently used scalars, vectors, and matrices in this paper} \label{notation}
\begin{center} 
\begin{tabular}{l|l} \toprule
$\ell$ & the maximal dimension of Krylov subspace generated with $\{\mbf A,\mbf r_0\}$\\ \noalign{\smallskip}\hline\noalign{\smallskip}
$m$ & the maximal dimension of Krylov subspace generated with $\{\mbf A,\mbf A\mbf r_0\}$\\ \noalign{\smallskip}\hline\noalign{\smallskip}
$\alpha=\ind(\mbf A)$ & the size of a largest Jordan block associated with zero eigenvalue of $\mbf A$\\ \noalign{\smallskip}\hline\noalign{\smallskip}
$\beta_1=\|\mbf r_0\|$ & the Euclidean norm of the initial residual vector $\mbf r_0$\\ \noalign{\smallskip}\hline\noalign{\smallskip}
$\wh\beta_1=\|\mbf A\mbf r_0\|$ & the Euclidean norm of the initial $\mbf A$-residual vector $\mbf A\mbf r_0$\\ \noalign{\smallskip}\hline\noalign{\smallskip}
$\kappa(\mbf A)=\|\mbf A\|\|\mbf A^\dag\|$ & the ratio of the largest singular value of $\mbf A$ to the smallest positive one\\ \noalign{\smallskip}\hline\noalign{\smallskip}
$\mbf A^\dag\mbf b$ & the pseudoinverse solution of $\mbf A\mbf x=\mbf b$\\ \noalign{\smallskip}\hline\noalign{\smallskip}
$\mbf A^\dag\mbf b+(\mbf I-\mbf A^\dag\mbf A)\mbf x_0$ & the orthogonal projection of $\mbf x_0$ onto the (least squares) solution set \\ \noalign{\smallskip}\hline\noalign{\smallskip}
$\mbf V_k$ $(1\leq k\leq \ell)$ & the matrix whose columns form an orthonormal basis of $\mcalk_k(\mbf A,\mbf r_0)$ \\ \noalign{\smallskip}\hline\noalign{\smallskip}
$\mbf H_{k+1,k}$ ($\mbf H_\ell$) & the matrix generated in the Arnoldi process for $\mcalk_k(\mbf A,\mbf r_0)$ \\ \noalign{\smallskip}\hline\noalign{\smallskip}
$\wh{\mbf V}_k$ $(1\leq k\leq m)$ & the matrix whose columns form an orthonormal basis of $\mcalk_k(\mbf A,\mbf A\mbf r_0)$\\ \noalign{\smallskip}\hline\noalign{\smallskip}
$\wh{\mbf H}_{k+1,k}$ ($\wh{\mbf H}_m$)  & the matrix generated in the Arnoldi process for $\mcalk_k(\mbf A,\mbf A\mbf r_0)$ \\ \noalign{\smallskip}
\bottomrule
\end{tabular}
\end{center}
\end{table}

\section{GMRES-type methods for singular range-symmetric linear systems}

\subsection{GMRES and a lifting strategy}
For any initial approximate solution $\mbf x_0$, at step $k$, GMRES determines the $k$th approximate solution \beq\label{gmressub}\mbf x_k:=\argmin_{\mbf x\in\mbf x_0+\mcalk_k(\mbf A,\mbf r_0)}\|\mbf b-\mbf A\mbf x\|.\eeq 
Since the columns of $\mbf V_k$ form an orthonormal basis of $\mcalk_k(\mbf A,\mbf r_0)$, using $\mbf A\mbf V_k=\mbf V_{k+1}\mbf H_{k+1,k}$ and $\mbf r_0=\beta_1\mbf v_1$, we have $\mbf x_k=\mbf x_0+\mbf V_k\mbf z_k$, where $\mbf z_k$ solves $\min_{\mbf z\in\mbbr^k}\|\beta_1\mbf e_1-\mbf H_{k+1,k}\mbf z\|.$ 

For singular $\mbf A$, Brown and Walker \cite{brown1997gmres} gave conditions under which the GMRES iterates converge safely to a least squares solution or to the pseudoinverse solution. More precisely, they proved the following results. (i) If $\ran(\mbf A)=\ran(\mbf A^\top)$ and $\mbf b\in\ran(\mbf A)$, then for all $0\leq k\leq\ell-1$, $\mbf x_k$ is not a solution, and $\mbf x_\ell=\mbf A^\dag\mbf b+(\mbf I-\mbf A^\dag\mbf A)\mbf x_0$, the orthogonal projection of $\mbf x_0$ onto the solution set $\{\mbf x\in\mbbr^n\ |\ \mbf A\mbf x=\mbf b\}$. (ii) If $\ran(\mbf A)=\ran(\mbf A^\top)$ and $\mbf b\notin\ran(\mbf A)$, then $\mbf x_{\ell-1}$ is a least squares solution of \eqref{lin}. 

Brown and Walker \cite{brown1997gmres} also studied the condition number of the upper-Hessenberg matrix $\mbf H_{k+1,k}$. They gave the following estimate. If $\ran(\mbf A)=\ran(\mbf A^\top)$ and $\mbf b\in\ran(\mbf A)$, then $\kappa(\mbf H_{k+1,k})\leq\kappa(\mbf A)$. Let $\mbf r_\star$ denote the least squares residual for \eqref{lin} and $\mbf r_k$ be the $k$th residual of GMRES. If $\ran(\mbf A)=\ran(\mbf A^\top)$ and $\mbf r_{k-1}\neq\mbf r_\star$, then \beq\label{kappa}\kappa(\mbf H_{k+1,k})\geq\frac{\|\mbf H_{k+1,k}\|}{\|\mbf A\|}\frac{\|\mbf r_{k-1}\|}{\sqrt{\|\mbf r_{k-1}\|^2-\|\mbf r_\star\|^2}}.\eeq The last estimate means that in the inconsistent range-symmetric case ($\mbf r_\star\neq\mbf 0$), the least squares problem \eqref{gmressub} becomes ill-conditioned as the GMRES iterate converges to a least squares solution.

Next we consider how to obtain the pseudoinvese solution for the case $\ran(\mbf A)=\ran(\mbf A^\top)$ and $\mbf b\notin\ran(\mbf A)$ from the final GMRES iterate $\mbf x_{\ell-1}$. Using the lifting strategy of \cite{liu2023obtai}, we define the lifted vector \beq\label{lifted} \wt{\mbf x}_{\ell-1}:=\mbf x_{\ell-1}-\frac{\mbf r_{\ell-1}^\top(\mbf x_{\ell-1}-\mbf x_0)}{\mbf r_{\ell-1}^\top\mbf r_{\ell-1}}\mbf r_{\ell-1},\eeq where $\mbf r_{\ell-1}:=\mbf b-\mbf A\mbf x_{\ell-1}$. We have the following result.

\begin{theorem}\label{main1}
If $\ran(\mbf A)=\ran(\mbf A^\top)$ and $\mbf b\notin\ran(\mbf A)$, then the lifted vector $\wt{\mbf x}_{\ell-1}$ in \eqref{lifted} is the orthogonal projection of $\mbf x_0$ onto the least squares solution set $\{\mbf x\in\mbbr^n\ |\ \mbf A^\top\mbf A\mbf x=\mbf A^\top\mbf b\}$. More precisely, we have $$\wt{\mbf x}_{\ell-1}=\mbf A^\dag\mbf b+(\mbf I-\mbf A^\dag\mbf A)\mbf x_0.$$ 
\end{theorem}
\proof 

It follows from $\mbf x_{\ell-1}$ is a least squares solution of \eqref{lin} that $\mbf r_{\ell-1}=(\mbf I-\mbf A\mbf A^\dag)\mbf b$ (see, for example, \cite[page 488]{meyer2023matri} for a proof). Since $\mbf x_{\ell-1}\in\mbf x_0+\mcalk_{\ell-1}(\mbf A,\mbf r_0)$, we can write $$\mbf x_{\ell-1}=\mbf x_0+\sum_{i=1}^{\ell-1}\alpha_i\mbf A^{i-1}\mbf r_0, \quad \alpha_i\in\mbbr.$$ 
Define $$\mbf f:=\alpha_1(\mbf I-\mbf A\mbf A^\dag)\mbf r_0=\alpha_1(\mbf I-\mbf A\mbf A^\dag)(\mbf b-\mbf A\mbf x_0)=\alpha_1(\mbf I-\mbf A\mbf A^\dag)\mbf b=\alpha_1\mbf r_{\ell-1}$$ 
and $$\mbf g:= \mbf x_{\ell-1}-\mbf f =\mbf x_0+\sum_{i=1}^{\ell-1}\alpha_i\mbf A^{i-1}\mbf r_0-\alpha_1(\mbf I-\mbf A\mbf A^\dag)\mbf r_0=\mbf x_0+\alpha_1\mbf A\mbf A^\dag\mbf r_0+ \sum_{i=2}^{\ell-1}\alpha_i\mbf A^{i-1}\mbf r_0.$$ 
The last two terms in last equation both lie in $\ran(\mbf A)$. Using $\mbf r_{\ell-1}\perp\ran(\mbf A)$, we get $\mbf f^\top\mbf g=\mbf f^\top(\mbf x_{\ell-1}-\mbf f)=\mbf f^\top\mbf x_0$. This gives $$\alpha_1\mbf r_{\ell-1}^\top\mbf x_{\ell-1}-\alpha_1^2\mbf r_{\ell-1}^\top\mbf r_{\ell-1}=\alpha_1\mbf r_{\ell-1}^\top\mbf x_0.$$ 
Since $\alpha_1=0$ implies $\mbf r_{\ell-1}\perp\mbf x_{\ell-1}-\mbf x_0$, we have $$\alpha_1=\mbf r_{\ell-1}^\top(\mbf x_{\ell-1}-\mbf x_0)/\mbf r_{\ell-1}^\top\mbf r_{\ell-1}.$$ 

Since $\ran(\mbf A)=\ran(\mbf A^\top)$, there exists a matrix $\mbf B\in\mbbr^{n\times n}$ satisfying $\mbf A^\top=\mbf A\mbf B$. Using $\mbf A^\top=\mbf A\mbf B$, $\mbf A\mbf A^\dag=(\mbf A\mbf A^\dag)^\top$, and $\mbf A\mbf A^\dag\mbf A=\mbf A$, we get $$\mbf A\mbf A\mbf A^\dag=((\mbf A\mbf A^\dag)^\top\mbf A^\top)^\top=(\mbf A\mbf A^\dag\mbf A\mbf B)^\top=(\mbf A\mbf B)^\top=\mbf A,$$ 
which implies $\mbf A\mbf f=\mbf 0.$ Thus we have $\mbf A\mbf g=\mbf A\mbf x_{\ell-1}$. This means that $\mbf A^\top\mbf A\mbf g=\mbf A^\top\mbf A\mbf x_{\ell-1}=\mbf A^\top\mbf b$, that is, $\mbf g$ is a least squares solution of \eqref{lin}. Now we write $$\mbf g=(\mbf I-\mbf A^\dag\mbf A)\mbf x_0+\mbf A^\dag\mbf A\mbf x_0+\alpha_1\mbf A\mbf A^\dag\mbf r_0+ \sum_{i=2}^{\ell-1}\alpha_i\mbf A^{i-1}\mbf r_0\in\mbf A^\dag\mbf b+\nul(\mbf A).$$ 
Since $\mbf A^\dag\mbf b\perp\nul(\mbf A)$, $(\mbf I-\mbf A^\dag\mbf A)\mbf x_0\in\nul(\mbf A)$, $\mbf A^\dag\mbf A\mbf x_0\in\ran(\mbf A^\top)$, and $\alpha_1\mbf A\mbf A^\dag\mbf r_0+ \sum_{i=2}^{\ell-1}\alpha_i\mbf A^{i-1}\mbf r_0\in\ran(\mbf A)=\ran(\mbf A^\top)$, we must have $$\mbf A^\dag\mbf A\mbf x_0+\alpha_1\mbf A\mbf A^\dag\mbf r_0+ \sum_{i=2}^{\ell-1}\alpha_i\mbf A^{i-1}\mbf r_0=\mbf A^\dag\mbf b,$$ 
which implies $$\wt{\mbf x}_{\ell-1}=\mbf x_{\ell-1}-\frac{\mbf r_{\ell-1}^\top(\mbf x_{\ell-1}-\mbf x_0)}{\mbf r_{\ell-1}^\top\mbf r_{\ell-1}}\mbf r_{\ell-1}=\mbf x_{\ell-1}-\alpha_1\mbf r_{\ell-1}=\mbf x_{\ell-1}-\mbf f=\mbf g=(\mbf I-\mbf A^\dag\mbf A)\mbf x_0+\mbf A^\dag\mbf b.$$ This completes the proof.
\endproof

    
\begin{corollary}\label{pseudo}
If $\ran(\mbf A)=\ran(\mbf A^\top)$, $\mbf b\notin\ran(\mbf A)$, and $\mbf x_0\in\ran(\mbf A)$, then the lifted vector $\wt{\mbf x}_{\ell-1}$ in \eqref{lifted} is the pseudoinverse solution $\mbf A^\dag\mbf b$.		
\end{corollary}
\proof Using $\mbf x_0\in\ran(\mbf A)=\ran(\mbf A^\top)$ and $(\mbf I-\mbf A^\dag\mbf A)\mbf A^\top=\mbf 0$, we have $\wt{\mbf x}_{\ell-1}=\mbf A^\dag\mbf b$. 
\endproof

Since the columns of $\mbf V_{\ell-1}$ form an orthonormal basis of $\mcalk_{\ell-1}(\mbf A,\mbf r_0)$, using $\mbf A\mbf V_{\ell-1}=\mbf V_\ell\mbf H_{\ell,\ell-1}$ and $\mbf r_0=\beta_1\mbf v_1$, we obtain $\mbf x_{\ell-1}=\mbf x_0+\mbf V_{\ell-1}\mbf z_{\ell-1}$, where $\mbf z_{\ell-1}$ solves $\min_{\mbf z\in\mbbr^{\ell-1}}\|\beta_1\mbf e_1-\mbf H_{\ell,\ell-1}\mbf z\|.$ If $\mbf A$ is skew-symmetric, i.e., $\mbf A^\top=-\mbf A$, then $\mbf H_\ell=\mbf V_\ell^\top\mbf A\mbf V_\ell$ is also skew-symmetric. The structure of $\mbf H_{\ell,\ell-1}$ yields that the odd entries of $\mbf z_{\ell-1}$ are zero (see, for example, \cite[section 8]{greif2016numer}). In this case, we have $$\mbf r_{\ell-1}^\top(\mbf x_{\ell-1}-\mbf x_0)=(\beta_1\mbf e_1-\mbf H_{\ell,\ell-1}\mbf z_{\ell-1})^\top\mbf V_\ell^\top\mbf V_{\ell-1}\mbf z_{\ell-1}=\beta_1\mbf e_1^\top\mbf z_{\ell-1}-\mbf z_{\ell-1}^\top\mbf H_\ell^\top\mbf z_{\ell-1}=0.$$ 
Hence, if $\mbf x_0\in\ran(\mbf A)$, $\mbf A^\top=-\mbf A$, and $\mbf b\notin\ran(\mbf A)$, then the $(\ell-1)$th GMRES iterate $\mbf x_{\ell-1}=\wt{\mbf x}_{\ell-1}=\mbf A^\dag\mbf b$. This result has been given in our previous work \cite [section 3.2]{du2023krylo}.


%

\subsection{RRGMRES}
 A variant of GMRES, named RRGMRES, was proposed in  \cite{calvetti2000gmres}. At step $k$, RRGMRES determines the $k$th approximate solution $$\mbf x_k^\rmr:=\argmin_{\mbf x\in\mbf x_0+\mbf \mcalk_k(\mbf A,\mbf A\mbf r_0)}\|\mbf b-\mbf A\mbf x\|^2.$$
Calvetti, Lewis, and Reichel \cite{calvetti2000gmres} proved that RRGMRES always determines the pseudoinverse solution if $\ran(\mbf A)=\ran(\mbf A^\top)$ and $\mbf x_0=\mbf 0$. More precisely, they proved the following results. (i) If $\mbf b\in\ran(\mbf A)$, $\ran(\mbf A)=\ran(\mbf A^\top)$, and $\mbf x_0=\mbf 0$, then $\mbf x_\ell^\rmr=\mbf A^\dag\mbf b$. (ii) If $\mbf b\notin\ran(\mbf A)$, $\ran(\mbf A)=\ran(\mbf A^\top)$, and $\mbf x_0=\mbf 0$, then $\mbf x_{\ell-1}^\rmr=\mbf A^\dag\mbf b$.
   
Since the columns of $\wh{\mbf V}_k$ form an orthonormal basis of $\mbf \mcalk_k(\mbf A,\mbf A\mbf r_0)$, using $\mbf A\wh{\mbf V}_k=\wh{\mbf V}_{k+1}\wh{\mbf H}_{k+1,k}$, we have \begin{align*}\min_{\mbf x\in\mbf x_0+\mbf \mcalk_k(\mbf A,\mbf A\mbf r_0)}\|\mbf b-\mbf A\mbf x\|^2&=\min_{\mbf z\in\mbbr^k}\|\mbf r_0-\mbf A\wh{\mbf V}_k\mbf z\|^2=\min_{\mbf z\in\mbbr^k}\|\mbf r_0-\wh{\mbf V}_{k+1}\wh{\mbf H}_{k+1,k}\mbf z\|^2\\ &=\min_{\mbf z\in\mbbr^k}\|\wh{\mbf V}_{k+1}^\top\mbf r_0-\wh{\mbf H}_{k+1,k}\mbf z\|^2+\|(\mbf I-\wh{\mbf V}_{k+1}\wh{\mbf V}_{k+1}^\top)\mbf r_0\|^2.\end{align*} 
Since $\mbf A\mbf r_0\in\ran(\mbf A)$, using the result of \cite{brown1997gmres}, we have $\kappa(\wh{\mbf H}_{k+1,k})\leq\kappa(\mbf A)$ if $\ran(\mbf A)=\ran(\mbf A^\top)$. Recall that the least squares problem \eqref{gmressub} of GMRES may become dangerously ill conditioned before a least squares is reached (see the estimate \eqref{kappa}). Therefore, for inconsistent range-symmetric linear systems, RRGMRES is a successful alternative to GMRES (see \cite{morikuni2018gmres} for examples and more discussion). 

\subsection{RSMAR: An iterative method for range-symmetric linear systems}
For range-symmetric linear systems, at step $k$ RSMAR generates an approximation $$\mbf x_k^\rma:=\argmin_{\mbf x\in\mbf x_0+\mcalk_k(\mbf A,\mbf r_0)}\|\mbf A(\mbf b-\mbf A\mbf x)\|.$$ 
Using $\mbf A\mbf V_k=\mbf V_{k+1}\mbf H_{k+1,k}$ and $\mbf r_0=\beta_1\mbf v_1$, we have 
\begin{align*}
	 \mbf A(\mbf b-\mbf A(\mbf x_0+\mbf V_k\mbf z))&=\mbf A\mbf r_0-\mbf A\mbf V_{k+1}\mbf H_{k+1,k}\mbf z \nn \\ & =\mbf V_{k+2}(\beta_1\mbf H_{k+2,k+1}\mbf e_1-\mbf H_{k+2,k+1}\mbf H_{k+1,k}\mbf z),\quad 1\leq k\leq \ell-2,\\
	 \mbf A(\mbf b-\mbf A(\mbf x_0+\mbf V_{\ell-1}\mbf z))&=\mbf A\mbf r_0-\mbf A\mbf V_\ell\mbf H_{\ell,\ell-1}\mbf z =\mbf V_\ell(\beta_1\mbf H_\ell\mbf e_1-\mbf H_\ell\mbf H_{\ell,\ell-1}\mbf z),\\
	\mbf A(\mbf b-\mbf A(\mbf x_0+\mbf V_{\ell}\mbf z))&=\mbf A\mbf r_0-\mbf A\mbf V_\ell\mbf H_\ell\mbf z =\mbf V_\ell(\beta_1\mbf H_\ell\mbf e_1-\mbf H_\ell^2\mbf z).
\end{align*} 	
Since the first $k$ columns of $\mbf V_\ell$ form an orthonormal basis of $\mcalk_k(\mbf A,\mbf r_0)$, we have $\mbf x_k^\rma=\mbf x_0+\mbf V_k\mbf z_k^\rma$, where $\mbf z_k^\rma$ solves the following subproblems of RSMAR \begin{subequations}
\begin{align}
	&\min_{\mbf z\in\mbbr^k}\|\beta_1\mbf H_{k+2,k+1}\mbf e_1-\mbf H_{k+2,k+1}\mbf H_{k+1,k}\mbf z\|,\quad 1\leq k\leq \ell-2, \label{suba}\\
	&\min_{\mbf z\in\mbbr^{\ell-1}}\|\beta_1\mbf H_\ell\mbf e_1-\mbf H_\ell\mbf H_{\ell,\ell-1}\mbf z\|,\\
	&\min_{\mbf z\in\mbbr^\ell}\|\beta_1\mbf H_\ell\mbf e_1-\mbf H_\ell^2\mbf z\|.
\end{align} 	
\end{subequations}
The following lemma is required to show that the RSMAR iterate $\mbf x_k^\rma$ for each $1\leq k\leq m$ (recall that $m$ given in \eqref{mm} is the maximal dimension of Krylov subspace generated with $\{\mbf A,\mbf A\mbf r_0\}$) is well defined.


\begin{lemma}\label{dimm} If $\ran(\mbf A)=\ran(\mbf A^\top)$ and $\mbf b\in\ran(\mbf A)$, then $\dim\mcalk_k(\mbf A,\mbf A^2\mbf r_0)=k$ for each $1\leq k\leq\ell$. If $\ran(\mbf A)=\ran(\mbf A^\top)$ and $\mbf b\notin\ran(\mbf A)$, then $\dim\mcalk_\ell(\mbf A,\mbf A^2\mbf r_0)=\ell-1$ and $\dim\mcalk_k(\mbf A,\mbf A^2\mbf r_0)=k$ for each $1\leq k\leq\ell-1$. 
\end{lemma}
\proof This is a direct result of \eqref{mm} and Lemma 3.1 of \cite{calvetti2000gmres}.
\endproof

Define $\mbf M_k=\mbf H_{k+2,k+1}\mbf H_{k+1,k}$ for $1\leq k\leq\ell-2$, $\mbf M_{\ell-1}=\mbf H_\ell\mbf H_{\ell,\ell-1}$, and $\mbf M_\ell=\mbf H_\ell^2$. Using Lemma \ref{dimm}, we next show that when $\mbf A$ is range-symmetric, $\mbf M_k$ for each $1\leq k\leq m$ has full column rank, which implies $\mbf z_k^\rma$ is unique for each $1\leq k\leq m$. We only consider the case $k=m=\ell-1$. All other cases are analogous. If $k=m=\ell-1$ (in this case $\mbf H_\ell$ is singular and $\mbf b\notin\ran(\mbf A)$), then $\dim\mcalk_{\ell-1}(\mbf A,\mbf A^2\mbf r_0)=\ell-1$ implies $\rank(\mbf A^2\mbf V_{\ell-1})=\rank(\mbf V_\ell\mbf H_\ell\mbf H_{\ell,\ell-1})=\rank(\mbf H_\ell\mbf H_{\ell,\ell-1})=\ell-1.$  Since $\mbf z_k^\rma$ is unique for each $1\leq k\leq m$, then $\mbf x_k^\rma=\mbf x_0+\mbf V_k\mbf z_k^\rma$ is well defined. Moreover, we have the following result.

\begin{theorem}\label{equiv} If $\ran(\mbf A)=\ran(\mbf A^\top)$ and $\mbf b\in\ran(\mbf A)$, then $\mbf x_\ell^\rma=\mbf x_\ell$. If $\ran(\mbf A)=\ran(\mbf A^\top)$ and $\mbf b\notin\ran(\mbf A)$, then $\mbf x_{\ell-1}^\rma=\mbf x_{\ell-1}$.
\end{theorem}
\proof When $\ran(\mbf A)=\ran(\mbf A^\top)$ and $\mbf b\in\ran(\mbf A)$, the matrix $\mbf H_\ell$ is invertible. So $\mbf z_\ell=\beta_1\mbf H_\ell^{-1}\mbf e_1=\mbf z_\ell^\rma$, which gives $\mbf x_\ell^\rma=\mbf x_\ell$. When $\ran(\mbf A)=\ran(\mbf A^\top)$ and $\mbf b\notin\ran(\mbf A)$, the matrix $\mbf H_\ell$ is singular and $\rank(\mbf H_\ell)=\ell-1$. It follows from $\ran(\mbf H_\ell\mbf H_{\ell,\ell-1})\subseteq\ran(\mbf H_\ell)$ and $\rank(\mbf H_\ell)=\rank(\mbf H_\ell\mbf H_{\ell,\ell-1})=\ell-1$ that $\ran(\mbf H_\ell\mbf H_{\ell,\ell-1})=\ran(\mbf H_\ell)$. This means that $\mbf H_\ell\mbf H_{\ell,\ell-1}\mbf z=\beta_1\mbf H_\ell\mbf e_1$ is consistent. Hence, we have $\mbf A(\mbf b-\mbf A\mbf x_{\ell-1}^\rma)=\mbf V_\ell(\beta_1\mbf H_\ell\mbf e_1-\mbf H_\ell\mbf H_{\ell,\ell-1}\mbf z_{\ell-1}^\rma)=\mbf 0$, which implies that $\mbf x_{\ell-1}^\rma\in\mbf x_0+\mcalk_{\ell-1}(\mbf A,\mbf r_0)$ is a least squares solution of \eqref{lin}. Since the final iterate GMRES iterate $\mbf x_{\ell-1}$ is also a least squares solution, by Theorem \ref{one}, it must hold that $\mbf x_{\ell-1}^\rma=\mbf x_{\ell-1}$. 
\endproof

Theorem \ref{equiv} means that for range-symmetric linear systems, GMRES and RSMAR terminate with the same least squares solution.

If $\ind(\mbf A)=1$, then the matrix $\wh{\mbf H}_m$ is invertible (see Theorem \ref{index1}). Hence, for each $1\leq k\leq m$, we have $\dim\mbf A^2(\mcalk_k(\mbf A, \mbf r_0))=\dim\mbf A(\mcalk_k(\mbf A, \mbf A\mbf r_0))=k$. Using similar analysis as before, we can conclude that the RSMAR iterate $\mbf x_k^\rma$ ($1\leq k\leq m$) is well defined when applied to linear systems with index one. Indeed, we have the following result.

\begin{theorem}\label{group}
If $\ind(\mbf A)=1$ and $\mbf b\in\ran(\mbf A)$, then $\mbf x_\ell^\rma=\mbf A^\rmd\mbf b+(\mbf I-\mbf A^\rmd\mbf A)\mbf x_0$. If $\ind(\mbf A)=1$ and $\mbf b\notin\ran(\mbf A)$, then $\mbf x_{\ell-1}^\rma$ satisfies $\mbf A^2\mbf x_{\ell-1}^\rma=\mbf A\mbf b$.
\end{theorem}
\proof Using $\mbf A\wh{\mbf V}_m=\wh{\mbf V}_m\wh{\mbf H}_m$ and $\mbf A\mbf r_0=\wh\beta_1\wh{\mbf v}_1$, we have 
\begin{align*}
	 \mbf A(\mbf b-\mbf A\mbf x_0-\wh{\mbf V}_m\wh{\mbf z})&=\mbf A\mbf r_0-\mbf A\wh{\mbf V}_m\wh{\mbf z} =\wh{\mbf V}_m(\wh\beta_1\mbf e_1-\wh{\mbf H}_m\wh{\mbf z}).
\end{align*} 	
Since the columns of $\wh{\mbf V}_m$ form an orthonormal basis of $\mcalk_m(\mbf A,\mbf A\mbf r_0)$, we have 
\begin{align*}
	\min_{\mbf x\in\mbf x_0+\mcalk_m(\mbf A,\mbf r_0)}\|\mbf A(\mbf b-\mbf A\mbf x)\|=&\min_{\wh{\mbf z}\in\mbbr^m}\|\wh\beta_1\mbf e_1-\wh{\mbf H}_m\wh{\mbf z}\|.
\end{align*} 	
When $\ind(\mbf A)=1$, the matrix $\wh{\mbf H}_m$ is invertible (see Theorem \ref{index1}). Then we have $$\min_{\mbf x\in\mbf x_0+\mcalk_m(\mbf A,\mbf r_0)}\|\mbf A(\mbf b-\mbf A\mbf x)\|=\min_{\wh{\mbf z}\in\mbbr^m}\|\wh\beta_1\mbf e_1-\wh{\mbf H}_m\wh{\mbf z}\|=0.$$ This means that $\mbf A(\mbf b-\mbf A\mbf x_m^\rma)=\mbf 0$. 

When $\ind(\mbf A)=1$ and $\mbf b\notin\ran(\mbf A)$, the matrix $\mbf H_\ell$ is singular and we have $m=\ell-1$. Therefore, $\mbf A^2\mbf x_{\ell-1}^\rma=\mbf A\mbf b$. 

When $\ind(\mbf A)=1$ and $\mbf b\in\ran(\mbf A)$, the matrix $\mbf H_\ell$ is nonsingular and we have $m=\ell$. So $\mbf b-\mbf A\mbf x_\ell^\rma\in\ran(\mbf A)\cap\nul(\mbf A)=\{\mbf 0\}$ (note that $\ran(\mbf A)\cap\nul(\mbf A)=\{\mbf 0\}$ is equivalent to $\ind(\mbf A)=1$). This means $\mbf x_\ell^\rma\in\mbf x_0+\mcalk_\ell(\mbf A,\mbf r_0)$ is a solution of $\mbf A\mbf x=\mbf b$. Using $\mbf r_0\in\ran(\mbf A)$, we have $\mbf b-\mbf A(\mbf x_0+\mbf A^\rmd\mbf r_0)=\mbf r_0-\mbf A\mbf A^\rmd\mbf r_0=\mbf 0$, which implies $\mbf x_0+\mbf A^\rmd\mbf r_0$ is also a solution of $\mbf A\mbf x=\mbf b$. By $\mbf x_0+\mbf A^\rmd\mbf r_0\in\mbf x_0+\mcalk_\ell(\mbf A,\mbf r_0)$ and Theorem \ref{one}, it must hold that $\mbf x_\ell^\rma=\mbf x_0+\mbf A^\rmd\mbf r_0=\mbf A^\rmd\mbf b+(\mbf I-\mbf A^\rmd\mbf A)\mbf x_0.$
\endproof

Since range-symmetric $\mbf A$ has index one and $\mbf A^2\mbf x=\mbf A\mbf b$ is equivalent to the normal equations $\mbf A^\top\mbf A\mbf x=\mbf A^\top\mbf b$ in the sense that they have the same solution set $\mbf A^\dag\mbf b+\nul(\mbf A)$, we know Theorem \ref{equiv} is a direct result of Theorems \ref{one} and \ref{group}.

Next, we provide two implementations for RSMAR, one based on the Arnoldi process for $\mcalk_k(\mbf A,\mbf A\mbf r_0)$ and the other based on the Arnoldi process for $\mcalk_k(\mbf A,\mbf r_0)$. 

\subsubsection{Implementation based on Arnoldi process for $\mcalk_k(\mbf A,\mbf A\mbf r_0)$}

The implementation discussed here is inspired by the approach proposed by Walker and Zhou \cite{walker1994simpl} for the implementation of GMRES. 

By $\mbf A\wh{\mbf V}_k=\wh{\mbf V}_{k+1}\wh{\mbf H}_{k+1,k}$ and $\mbf A\mbf r_0=\wh\beta_1\wh{\mbf v}_1$, we have $$\mbf A(\mbf b-\mbf A\mbf x_0-\wh{\mbf V}_k\wh{\mbf z})=\mbf A\mbf r_0-\mbf A\wh{\mbf V}_k\wh{\mbf z} =\wh{\mbf V}_{k+1}(\wh\beta_1\mbf e_1-\wh{\mbf H}_{k+1,k}\wh{\mbf z}).$$
Since the columns of $\wh{\mbf V}_k$ form an orthonormal basis of $\mcalk_k(\mbf A,\mbf A\mbf r_0)$, we have \beq\label{sub1}\min_{\mbf x\in\mbf x_0+\mcalk_k(\mbf A,\mbf r_0)}\|\mbf A(\mbf b-\mbf A\mbf x)\|=\min_{\wh{\mbf z}\in\mbbr^k}\|\wh\beta_1\mbf e_1-\wh{\mbf H}_{k+1,k}\wh{\mbf z}\|.\eeq Now we introduce the QR factorization $$\wh{\mbf H}_{k+1,k}=\wh{\mbf Q}_{k+1}\bem\wh{\mbf R}_k\\ \mbf 0\eem,$$ 
where $\wh{\mbf Q}_{k+1}\in\mbbr^{(k+1)\times(k+1)}$ is orthogonal and upper Hessenberg, and $\wh{\mbf R}_k\in\mbbr^{k\times k}$ is nonsingular and upper triangular. Define $\wh{\mbf t}_{k+1}:=\wh{\mbf Q}_{k+1}^\top\wh\beta_1\mbf e_1$. The vector $\wh{\mbf z}_k:=\wh{\mbf R}_k^{-1}\bem\mbf I_k & \mbf 0\eem\wh{\mbf t}_{k+1}$ solves the least squares problem in the right hand side of \eqref{sub1}. Note that $\mcalk_k(\mbf A,\mbf r_0)=\spa\{\mbf r_0,\wh{\mbf v}_1,\ldots,\wh{\mbf v}_{k-1}\}$. The RSMAR iterate $\mbf x_k^\rma$ can be expressed as $$\mbf x_k^\rma=\mbf x_0+\bem \mbf r_0 &\wh{\mbf V}_{k-1}\eem \mbf z_k,$$ 
where $\mbf z_k$ solves $$\mbf A\bem \mbf r_0 & \wh{\mbf V}_{k-1}\eem\mbf z=\bem \wh\beta_1\wh{\mbf V}_k\mbf e_1 &\wh{\mbf V}_k\wh{\mbf H}_{k,k-1}\eem\mbf z=\wh{\mbf V}_k\bem \wh\beta_1\mbf e_1 &\wh{\mbf H}_{k,k-1}\eem\mbf z=\wh{\mbf V}_k\wh{\mbf z}_k.$$ 
Define $\wt{\mbf R}_k:=\bem \wh\beta_1\mbf e_1 &\wh{\mbf H}_{k,k-1}\eem$, which is upper triangular and invertible. We finally have \beq\label{imple1}\mbf x_k^\rma=\mbf x_0+\bem \mbf r_0 &\wh{\mbf V}_{k-1}\eem\wt{\mbf R}_k^{-1}\wh{\mbf z}_k.\eeq 
Note that we also have $$\|\mbf A\mbf r_k^\rma\|=\|\mbf A(\mbf b-\mbf A\mbf x_k^\rma)\|=\|\wh\beta_1\mbf e_1-\wh{\mbf H}_{k+1,k}\wh{\mbf z}_k\|=|\mbf e_{k+1}^\top\wh{\mbf t}_{k+1}|.$$ 
The approach given above is summarized as Algorithm 1. 

\begin{table}[H]
\centering
\begin{tabular*}{170mm}{l}
\toprule {\bf Algorithm 1}. RSMAR-I: implementation based on $\mbf A\wh{\mbf V}_k=\wh{\mbf V}_{k+1}\wh{\mbf H}_{k+1,k}$ for $\mcalk_k(\mbf A,\mbf A\mbf r_0)$
\\ \hline\noalign{\smallskip} \quad {\bf Require}: $\mbf A\in\mbbr^{n\times n}$ with $\ran(\mbf A)=\ran(\mbf A^\top)$, $\mbf b\in\mbbr^n$, $\mbf x_0\in\mbbr^n$, ${\tt tol}>0$, ${\tt maxit}>0$\hspace{10mm}
\\ \noalign{\smallskip} \quad\hspace{.64mm} 1:\ \ $\mbf r_0:=\mbf b-\mbf A\mbf x_0$, $\wh\beta_1:=\|\mbf A\mbf r_0\|$. If $\wh\beta_1<\tt tol$, accept $\mbf x_0$ and exit. 
\\ \noalign{\smallskip} \quad\hspace{.64mm} 2:\ \  $\wh{\mbf v}_1:=\mbf A\mbf r_0/\wh\beta_1$
\\ \noalign{\smallskip} \quad\hspace{.64mm} 3:\ \  {\bf for} $k=1, 2, \ldots, \tt maxit$ {\bf do}
\\ \noalign{\smallskip} \quad\hspace{.64mm} 4:\ \
 \qquad $\wh{\mbf v}_{k+1}:=\mbf A\wh{\mbf v}_k$ 
\\ \noalign{\smallskip} \quad\hspace{.64mm} 5:\ \ \qquad {\bf for} $i=1,2,\ldots,k$ {\bf do} 
\\ \noalign{\smallskip}\hspace{.64mm} \quad 6:\ \ \qquad\qquad $\wh h_{ik}:=\wh{\mbf v}_i^\top\wh{\mbf v}_{k+1}$  
\\ \noalign{\smallskip} \quad\hspace{.64mm} 7:\ \ \qquad\qquad $\wh{\mbf v}_{k+1}:=\wh{\mbf v}_{k+1}-\wh h_{ik}\wh{\mbf v}_i$ 
\\ \noalign{\smallskip} \quad\hspace{.64mm} 8:\ \ \qquad {\bf end}
\\ \noalign{\smallskip} \quad\hspace{.64mm} 9:\ \ \qquad $\wh h_{k+1,k}:=\|\wh{\mbf v}_{k+1}\|$ 
\\ \noalign{\smallskip} \quad 10:\ \ \qquad $\wh{\mbf v}_{k+1}:=\wh{\mbf v}_{k+1}/\wh h_{k+1,k}$
\\ \noalign{\smallskip} \quad 11:\ \  \qquad $\wh{\mbf H}_{k+1,k}:=\bem \wh{\mbf H}_{k,k-1} & \wh{\mbf h}_k\\\mbf 0 & \wh h_{k+1,k} \eem$ with $\wh{\mbf h}_k:=\bem \wh h_{1k} &\cdots & \wh h_{kk}\eem^\top$
\\ \noalign{\smallskip} \quad 12:\ \  \qquad $\wh{\mbf Q}_{k+1}\bem\wh{\mbf R}_k\\ \mbf 0\eem=\wh{\mbf H}_{k+1,k}$ \hfill  {\color[gray]{0.5} QR factorization of $\wh{\mbf H}_{k+1,k}$}
\\ \noalign{\smallskip} \quad 13:\ \  \qquad $\wh{\mbf t}_{k+1}:=\wh{\mbf Q}_{k+1}^\top\wh\beta_1\mbf e_1$
\\ \noalign{\smallskip} \quad 14:\ \  \qquad $\rho_k:=|\mbf e_{k+1}^\top\wh{\mbf t}_{k+1}|$
\\ \noalign{\smallskip} \quad 15:\ \  \qquad {\bf if} $\rho_k<\tt tol$ {\bf then}
\\ \noalign{\smallskip} \quad 16:\ \  \qquad\qquad $\wh{\mbf z}_k:=\wh{\mbf R}_k^{-1}\bem\mbf I_k & \mbf 0\eem\wh{\mbf t}_{k+1}$
\\ \noalign{\smallskip} \quad 17:\ \  \qquad\qquad $\wt{\mbf R}_k:=\bem \wh\beta_1\mbf e_1 &\wh{\mbf H}_{k,k-1}\eem $
\\ \noalign{\smallskip} \quad 18:\ \  \qquad\qquad $\mbf x_k^\rma:=\mbf x_0+\bem\mbf r_0 & \wh{\mbf V}_{k-1}\eem\wt{\mbf R}_k^{-1}\wh{\mbf z}_k$ 
\\ \noalign{\smallskip} \quad 19:\ \  \qquad\qquad Accept $\mbf x_k$ and exit.
\\ \noalign{\smallskip} \quad 20:\ \  \qquad {\bf end if}
\\ \noalign{\smallskip} \quad 21:\ \ {\bf end for} 
\\ 
\bottomrule
\end{tabular*}
\end{table}




\subsubsection{Implementation based on Arnoldi process for $\mcalk_k(\mbf A,\mbf r_0)$}

The implementation discussed here is inspired by the approach proposed by Neuman, Reichel, and Sadok  \cite{neuman2012algor,neuman2012imple} for the implementation of RRGMRES. 

We first introduce the QR factorization $$\mbf H_{k+1,k}=\mbf Q_{k+1}\bem\mbf R_k\\ \mbf 0\eem,$$ where $\mbf Q_{k+1}\in\mbbr^{(k+1)\times(k+1)}$ is orthogonal and upper Hessenberg, and $\mbf R_k\in\mbbr^{k\times k}$ is nonsingular and upper triangular. The subproblem \eqref{suba} of RSMAR can be written as \begin{align}\label{subat}\min_{\mbf z\in\mbbr^k}\|\beta_1\mbf H_{k+2,k+1}\mbf e_1-\mbf H_{k+2,k+1}\mbf H_{k+1,k}\mbf z\|=\min_{\wt{\mbf z}\in\mbbr^k}\|\beta_1(h_{11}\mbf e_1+h_{21}\mbf e_2)-\mbf H_{k+2,k+1}\mbf Q_{k+1}\bem\mbf I_k\\ \mbf 0\eem\wt{\mbf z}\|.\end{align}
The matrix $$\wt{\mbf H}_{k+2,k}:=\mbf H_{k+2,k+1}\mbf Q_{k+1}\bem\mbf I_k\\ \mbf 0\eem\in\mbbr^{(k+2)\times k}$$ vanishes below the sub-subdiagonal because $\mbf H_{k+2,k+1}$ and $\mbf Q_{k+1}$ are both upper Hessenberg. We then introduce the QR factorization $$\wt{\mbf H}_{k+2,k}=\wt{\mbf Q}_{k+2}\bem\wt{\mbf R}_k\\ \mbf 0\eem,$$ where $\wt{\mbf Q}_{k+2}\in\mbbr^{(k+2)\times(k+2)}$ is orthogonal and $\wt{\mbf R}_k\in\mbbr^{k\times k}$ is nonsingular and upper triangular. 
Define $\wt{\mbf t}_{k+2}:=\wt{\mbf Q}_{k+2}^\top\beta_1(h_{11}\mbf e_1+h_{21}\mbf e_2)\in\mbbr^{k+2}$. The vector $\wt{\mbf z}_k:=\wt{\mbf R}_k^{-1}\bem\mbf I_k & \mbf 0\eem\wt{\mbf t}_{k+2}$ solves the least squares problem in the right hand side of \eqref{subat}, and the vector $\mbf z_k:=\mbf R_k^{-1}\wt{\mbf z}_k$ solves the least squares problem in the left hand side of \eqref{subat}. Hence the RSMAR iterate $\mbf x_k^\rma$ can be expressed as $$\mbf x_k^\rma=\mbf x_0+\mbf V_k\mbf z_k=\mbf x_0+\mbf V_k\mbf R_k^{-1}\wt{\mbf z}_k.$$ 
Note that we also have $$\|\mbf A\mbf r_k^\rma\|=\|\mbf A(\mbf b-\mbf A\mbf x_k^\rma)\|=\|\beta_1\mbf H_{k+2,k+1}\mbf e_1-\mbf H_{k+2,k+1}\mbf H_{k+1,k}\mbf z_k\|=\sqrt{(\mbf e_{k+1}^\top\wt{\mbf t}_{k+2})^2+(\mbf e_{k+2}^\top\wt{\mbf t}_{k+2})^2}.$$ 
The approach given above is summarized as Algorithm 2. 

\begin{table}[H]
\centering
\begin{tabular*}{170mm}{l}
\toprule {\bf Algorithm 2}. RSMAR-II: implementation based on $\mbf A \mbf V_k= \mbf V_{k+1}\mbf H_{k+1,k}$ for $\mcalk_k(\mbf A,\mbf r_0)$
\\ \hline\noalign{\smallskip} \quad {\bf Require}: $\mbf A\in\mbbr^{n\times n}$ with $\ran(\mbf A)=\ran(\mbf A^\top)$, $\mbf b\in\mbbr^n$, $\mbf x_0\in\mbbr^n$, ${\tt tol}>0$, ${\tt maxit}>0$ \hspace{10mm}
\\ \noalign{\smallskip} \quad\hspace{.64mm} 1:\ \ $\mbf r_0:=\mbf b-\mbf A\mbf x_0$, $\beta_1:=\|\mbf r_0\|$, $\wh\beta_1:=\|\mbf A\mbf r_0\|$. If $\wh\beta_1<\tt tol$, accept $\mbf x_0$ and exit. 
\\ \noalign{\smallskip} \quad\hspace{.64mm} 2:\ \  $\mbf v_1:=\mbf r_0/\beta_1$, 
\\ \noalign{\smallskip} \quad\hspace{.64mm} 3:\ \  $\mbf v_2:=\mbf A\mbf v_1$,  $h_{11}:=\mbf v_1^\top\mbf v_2$, $\mbf v_2:=\mbf v_2-h_{11}\mbf v_1$, $h_{21}:=\|\mbf v_2\|$, $\mbf v_2:=\mbf v_2/h_{21}$, $\mbf H_{2,1}:=\bem h_{11}\\ h_{21}\eem$
\\ \noalign{\smallskip} \quad\hspace{.64mm} 4:\ \  {\bf for} $k=1, 2, \ldots, \tt maxit$ {\bf do}
\\ \noalign{\smallskip} \quad\hspace{.64mm} 5:\ \
 \qquad $\mbf v_{k+2}:=\mbf A\mbf v_{k+1}$ 
\\ \noalign{\smallskip} \quad\hspace{.64mm} 6:\ \ \qquad {\bf for} $i=1,2,\ldots,k+1$ {\bf do} 
\\ \noalign{\smallskip}\hspace{.64mm} \quad 7:\ \ \qquad\qquad $h_{i,k+1}:=\mbf v_i^\top\mbf v_{k+2}$  
\\ \noalign{\smallskip} \quad\hspace{.64mm} 8:\ \ \qquad\qquad $\mbf v_{k+2}:=\mbf v_{k+2}-h_{i,k+1}\mbf v_i$ 
\\ \noalign{\smallskip} \quad\hspace{.64mm} 9:\ \ \qquad {\bf end}
\\ \noalign{\smallskip} \quad 10:\ \ \qquad $h_{k+2,k+1}:=\|\mbf v_{k+2}\|$ 
\\ \noalign{\smallskip} \quad 11:\ \ \qquad $\mbf v_{k+2}:=\mbf v_{k+2}/h_{k+2,k+1}$
\\ \noalign{\smallskip} \quad 12:\ \  \qquad $\mbf H_{k+2,k+1}:=\bem \mbf H_{k+1,k} & \mbf h_{k+1}\\\mbf 0 & h_{k+2,k+1} \eem$ with $\mbf h_{k+1}:=\bem h_{1,k+1} &\cdots & h_{k+1,k+1}\eem^\top$
\\ \noalign{\smallskip} \quad 13:\ \  \qquad $\mbf Q_{k+1}\bem\mbf R_k\\ \mbf 0\eem=\mbf H_{k+1,k}$ \hfill  {\color[gray]{0.5} QR factorization of $\mbf H_{k+1,k}$}
\\ \noalign{\smallskip} \quad 14:\ \  \qquad $\wt{\mbf Q}_{k+2}\bem\wt{\mbf R}_k\\ \mbf 0\eem=\mbf H_{k+2,k+1}\mbf Q_{k+1}\bem\mbf I_k\\ \mbf 0\eem$ \hfill  {\color[gray]{0.5} QR factorization of $\mbf H_{k+2,k+1}\mbf Q_{k+1}\bem\mbf I_k\\ \mbf 0\eem$}
\\ \noalign{\smallskip} \quad 15:\ \ \qquad $\wt{\mbf t}_{k+2}:=\wt{\mbf Q}_{k+2}^\top\beta_1(h_{11}\mbf e_1+h_{21}\mbf e_2)$
\\ \noalign{\smallskip} \quad 16:\ \ \qquad $\rho_k=\sqrt{(\mbf e_{k+1}^\top\wt{\mbf t}_{k+2})^2+(\mbf e_{k+2}^\top\wt{\mbf t}_{k+2})^2}$
\\ \noalign{\smallskip} \quad 17:\ \  \qquad {\bf if} $\rho_k<\tt tol$ {\bf then}
\\ \noalign{\smallskip} \quad 18:\ \  \qquad\qquad $\wt{\mbf z}_k:=\wt{\mbf R}_k^{-1}\bem\mbf I_k & \mbf 0\eem\wt{\mbf t}_{k+2}$
\\ \noalign{\smallskip} \quad 19:\ \  \qquad\qquad $\mbf x_k:=\mbf x_0+\mbf V_k\mbf R_k^{-1}\wt{\mbf z}_k$ 
\\ \noalign{\smallskip} \quad 20:\ \  \qquad\qquad Accept $\mbf x_k$ and exit.
\\ \noalign{\smallskip} \quad 21:\ \  \qquad {\bf end if}
\\ \noalign{\smallskip} \quad 22:\ \ {\bf end for} 
\\ 
\bottomrule
\end{tabular*}
\end{table}

\subsection{DGMRES}
DGMRES is a GMRES-type method for the Drazin-inverse solution of consistent and inconsistent linear systems $\mbf A\mbf x=\mbf b$. At step $k$, DGMRES determines the $k$th approximate solution $$\mbf x_k^\rmd:=\argmin_{\mbf x\in\mbf x_0+\mcalk_k(\mbf A,\mbf A^\alpha\mbf r_0)}\|\mbf A^\alpha(\mbf b-\mbf A\mbf x)\|,$$ where $\alpha={\rm ind}(\mbf A)$ is the index of $\mbf A$. Sidi \cite{sidi1999unifi,sidi2001dgmre} proved that DGMRES always determines the Drazin-inverse solution if $\mbf x_0=\mbf 0$. More precisely, we have the following results. (i) If $\ind(\mbf A)=1$, $\mbf b\in\ran(\mbf A)$, and $\mbf x_0=\mbf 0$, then $\mbf x_\ell^\rmd=\mbf A^\rmd\mbf b$. (ii) If $\ind(\mbf A)=1$, $\mbf b\notin\ran(\mbf A)$, and $\mbf x_0=\mbf 0$, then $\mbf x_{\ell-1}^\rmd=\mbf A^\rmd\mbf b$.

Since range-symmetric $\mbf A$ has index one and satisfies $\mbf A^\dag=\mbf A^\rmd$, DGMRES applied to range-symmetric linear systems always determines the pseudoinverse solution. Actually, DGMRES applied to range-symmetric linear systems can be viewed as a range restricted RSMAR method since the minimization problem is 
$$\min_{\mbf x\in\mbf x_0+\mcalk_k(\mbf A,\mbf A\mbf r_0)}\|\mbf A(\mbf b-\mbf A\mbf x)\|.$$


\subsection{Summary of GMRES-type methods for the pseudoinvese solution}

We summarize the four methods (GMRES, RRGMRES, RSMAR, and DGMRES) discussed in this section in Table \ref{sg}. We use $\mbf x_0=\mbf 0$ and focus on their final iterate when applied to range-symmetric linear systems. Both consistent and in consistent cases are included. We have the following results. \bit
\item For the consistent case, the four methods terminate at step $\ell$, and give the pseudoinverse solution.
\item For the inconsistent case, the four methods terminate at step $\ell-1$. RRGMRES and DGMRES give the pseudoinverse solution. GMRES and RSMAR terminate with the same least squares solution (see Theorem \ref{equiv}). The lifting strategy \eqref{lifted} can be used to get the pseudoinverse solution. 
\item GMRES and RRGMRES have residual minimization property and the residual norm is nonincreasing. 
RSMAR and DGMRES have $\mbf A$-residual minimization property and the $\mbf A$-residual norm is nonincreasing. 
\eit 

\begin{table}[H]
\caption{Minimization property and final iterate of GMRES-type methods for the pseudoinvese solution of range-symmetric linear systems.} \label{sg}
\begin{center} 
\begin{tabular}{c|c|c|c} \toprule
Method  & Minimization property at step $k$ & Consistent case &  Inconsistent case\\ 
\noalign{\smallskip}\hline\noalign{\smallskip}
GMRES & $\mbf x_k=\argmin_{\mbf x\in\mcalk_k(\mbf A,\mbf b)}\|\mbf b-\mbf A\mbf x\|$ & $\mbf x_\ell=\mbf A^\dag\mbf b$ & $\mbf A\mbf r_{\ell-1}=\mbf 0$, $\wt{\mbf x}_{\ell-1}=\mbf A^\dag\mbf b$ \\ \noalign{\smallskip}\hline\noalign{\smallskip}
RRGMRES & $\mbf x_k^\rmr=\argmin_{\mbf x\in\mcalk_k(\mbf A,\mbf A\mbf b)}\|\mbf b-\mbf A\mbf x\|$  & $\mbf x_\ell^\rmr=\mbf A^\dag\mbf b$ & $\mbf x_{\ell-1}^\rmr=\mbf A^\dag\mbf b$ \\ \noalign{\smallskip}\hline\noalign{\smallskip}
RSMAR & $\mbf x_k^\rma=\argmin_{\mbf x\in\mcalk_k(\mbf A,\mbf b)}\|\mbf A(\mbf b-\mbf A\mbf x)\|$  & $\mbf x_\ell^\rma=\mbf A^\dag\mbf b$ & $\mbf A\mbf r_{\ell-1}^\rma=\mbf 0$, $\wt{\mbf x}_{\ell-1}^\rma=\mbf A^\dag\mbf b$ \\ \noalign{\smallskip}\hline\noalign{\smallskip}
DGMRES & $\mbf x_k^\rmd=\argmin_{\mbf x\in\mcalk_k(\mbf A,\mbf A\mbf b)}\|\mbf A(\mbf b-\mbf A\mbf x)\|$ & $\mbf x_\ell^\rmd=\mbf A^\dag\mbf b$ & $\mbf x_{\ell-1}^\rmd=\mbf A^\dag\mbf b$ \\ \noalign{\smallskip}\bottomrule
\end{tabular}
\end{center}
\end{table}

\section{MINRES-type methods for singular symmetric linear systems}
In this section, we assume that $\mbf A$ is symmetric, i.e., $\mbf A^\top=\mbf A$. The matrix $\mbf H_\ell$ in \eqref{ell} is symmetric and tridiagonal, and it is nonsingular if and only if $\mbf b\in\ran(\mbf A)$ \cite[section 2.1 property 4]{choi2011minre}. For simplicity, in the following discussion, we choose $\mbf x_0=\mbf 0$. GMRES applied to symmetric linear systems is theoretically equivalent to MINRES \cite{paige1975solut}, which has short recurrences. 

\subsection{MINRES and a lifting strategy}
The subproblems of MINRES are 
 \begin{subequations} 
 \begin{align*} \min_{\mbf x\in\mcalk_k(\mbf A,\mbf b)}\|\mbf b-\mbf A\mbf x\|&=\min_{\mbf z\in\mbbr^k}\|\beta_1\mbf e_1-\mbf H_{k+1,k}\mbf z\|, \quad 1\leq k\leq\ell-1,\\ \min_{\mbf x\in\mcalk_\ell(\mbf A,\mbf b)}\|\mbf b-\mbf A\mbf x\|&=\min_{\mbf z\in\mbbr^\ell}\|\beta_1\mbf e_1-\mbf H_\ell\mbf z\|.
\end{align*}\end{subequations}
At step $k$, MINRES minimizes $\|\mbf r_k\|$ over $\mcalk_k(\mbf A,\mbf r_0)$ but not $\|\mbf A\mbf r_k\|$. If $\mbf b\in\ran(\mbf A)$, then the $\ell$th MINRES iterate $\mbf x_\ell$ is the pseudoinverse solution (see \cite[Theorem 3.1]{choi2011minre}). If $\mbf b\notin\ran(\mbf A)$, then $\mbf H_\ell$ is singular with $\rank(\mbf H_\ell)=\ell-1$, and the $(\ell-1)$th MINRES iterate $\mbf x_{\ell-1}$ is a least squares solution, but not necessarily the pseudoinverse solution (see \cite[Theorem 3.2]{choi2011minre}). Liu, Milzarek, and Roosta \cite[Theorem 1]{liu2023obtai} proved that the lifted vector $$\wt{\mbf x}_{\ell-1}=\mbf x_{\ell-1}-\frac{\mbf r_{\ell-1}^\top\mbf x_{\ell-1}}{\mbf r_{\ell-1}^\top\mbf r_{\ell-1}}\mbf r_{\ell-1}$$ is the pseudoinverse solution. 

%
%

\subsection{MINARES}
The $k$th iterate of MINARES, denoted by $\mbf x_k^\rma$, solves $$\min_{\mbf x\in\mcalk_k(\mbf A,\mbf b)}\|\mbf A(\mbf b-\mbf A\mbf x)\|=\min_{\mbf x\in\mcalk_k(\mbf A,\mbf b)}\|\mbf A\mbf b-\mbf A^2\mbf x\|.$$ If $\mbf b\in\ran(\mbf A)$, then the $\ell$th MINARES iterate $\mbf x_\ell^\rma$ is the pseudoinverse solution (see \cite[Theorem 4.4]{montoison2023minar}). Hence in this case, $\mbf x_\ell^\rma$ coincides with the $\ell$th MINRES iterate $\mbf x_\ell$. If $\mbf b\notin\ran(\mbf A)$, then the $(\ell-1)$th MINARES iterate $\mbf x_{\ell-1}^\rma$ is a least squares solution (see \cite[Theorem 4.5]{montoison2023minar}). The following theorem is a direct result of Theorem \ref{equiv} because RSMAR and MINARES are theoretically equivalent for symmetric linear systems. Here, we would like to provide a direct proof rather than using Theorem \ref{equiv}.

\begin{theorem}\label{main2}
If $\mbf A^\top=\mbf A$ and $\mbf b\notin\ran(\mbf A)$, then the $(\ell-1)$th MINARES iterate $\mbf x_{\ell-1}^\rma$ is equal to the $(\ell-1)$th MINRES iterate $\mbf x_{\ell-1}$.
\end{theorem}


\proof Using $\mbf A\mbf V_{\ell-1}=\mbf V_\ell\mbf H_{\ell,\ell-1}$, $\mbf b=\beta_1\mbf v_1$, and $\mbf V_\ell^\top\mbf V_\ell=\mbf I$, we obtain $\mbf x_{\ell-1}=\mbf V_{\ell-1}\mbf z_{\ell-1}$, where $\mbf z_{\ell-1}$ solves $$\min_{\mbf z\in\mbbr^{\ell-1}}\|\beta_1\mbf e_1-\mbf H_{\ell,\ell-1}\mbf z\|.$$
Similarly, we have $\mbf x_{\ell-1}^\rma=\mbf V_{\ell-1}\mbf z_{\ell-1}^\rma$, where $\mbf z_{\ell-1}^\rma$ solves $$\min_{\mbf z\in\mbbr^{\ell-1}}\|\beta_1\mbf H_\ell\mbf e_1-\mbf H_\ell\mbf H_{\ell,\ell-1}\mbf z\|.$$ 

Next we show that $\mbf z_{\ell-1}^\rma=\mbf z_{\ell-1}$, which yields $\mbf x_{\ell-1}^\rma=\mbf x_{\ell-1}$. Since $\rank(\mbf H_\ell)=\ell-1$ and $\mbf H_{\ell}$ is symmetric, then $\mbf H_\ell$ has a decomposition  $\mbf H_\ell=\mbf U_{\ell-1}\mbf\Lambda_{\ell-1}\mbf U_{\ell-1}^\top$, where $\mbf \Lambda_{\ell-1}$ is a diagonal matrix with nonzero eigenvalues of $\mbf H_\ell$ as diagonal entries, and $\mbf U_{\ell-1}$ is an $\ell\times(\ell-1)$ matrix with corresponding unit eigenvectors of $\mbf H_\ell$ as columns. It follows from $\ran(\mbf U_{\ell-1})=\ran(\mbf H_\ell)=\ran(\mbf H_{\ell,\ell-1})$ that there exists a nonsingular matrix $\mbf C_{\ell-1}\in\mbbr^{(\ell-1)\times(\ell-1)}$ such that $\mbf H_{\ell,\ell-1}=\mbf U_{\ell-1}\mbf C_{\ell-1}$. Then it follows $$\mbf z_{\ell-1}=\beta_1\mbf C_{\ell-1}^{-1}\mbf U_{\ell-1}^\top\mbf e_1=\beta_1(\mbf H_\ell\mbf H_{\ell,\ell-1})^\dag\mbf H_\ell\mbf e_1=\mbf z_{\ell-1}^\rma.$$ This completes the proof.
\endproof

Theorem \ref{main2} means that for the case $\mbf A=\mbf A^\top$ and $\mbf b\notin\ran(\mbf A)$, MINARES and MINRES terminate with the same least squares solution. The MINARES implementation (only short recurrences are required) in \cite[section 4]{montoison2023minar} determines $\mbf x_{\ell-1}^\rma$ in exact arithmetic, and can not find the pseudoinverse solution. By Theorem \ref{main2}, the lifted vector $$\wt{\mbf x}_{\ell-1}^\rma=\mbf x_{\ell-1}^\rma-\frac{(\mbf r_{\ell-1}^\rma)^\top\mbf x_{\ell-1}^\rma}{(\mbf r_{\ell-1}^\rma)^\top\mbf r_{\ell-1}^\rma}\mbf r_{\ell-1}^\rma$$ is the pseudoinverse solution. Here, $\mbf r_{\ell-1}^\rma:=\mbf b-\mbf A \mbf x_{\ell-1}^\rma.$

\subsubsection{A new implementation of MINARES}
MINARES is mathematically equivalent to RSMAR applied to symmetric linear systems. The MINARES implementation in \cite[section 4]{montoison2023minar} is based on the Arnoldi relation $\mbf A\mbf V_k=\mbf V_{k+1}\mbf H_{k+1,k}$, and thus can be viewed as a short recurrence variant of RSMAR-II (Algorithm 2). Now we derive a new implementation of MINARES, which is based on $\mbf A\wh{\mbf V}_k=\wh{\mbf V}_{k+1}\wh{\mbf H}_{k+1,k}$ and can be viewed as a short recurrence variant of RSMAR-I (Algorithm 1).

If $\mbf A=\mbf A^\top$, then the matrix $\wh{\mbf H}_k$ is symmetric and tridiagonal. The Arnoldi process reduces to the Lanczos process \cite{lanczos1950itera}. After $k$ iterations, we have  
 \begin{align*} \mbf A\wh{\mbf V}_k=\wh{\mbf V}_{k+1}\wh{\mbf H}_{k+1,k}=\wh{\mbf V}_k\wh{\mbf H}_k+\wh\beta_{k+1}\wh{\mbf v}_{k+1}\mbf e_k^\top,
\end{align*}
where $$\wh{\mbf V}_k=\bem\wh{\mbf v}_1 & \wh{\mbf v}_2 & \cdots & \wh{\mbf v}_k\eem,\qquad \wh{\mbf H}_k=\bem \wh \alpha_1 & \wh \beta_2 & & \\ \wh \beta_2 & \wh \alpha_2 & \ddots &  \\ & \ddots & \ddots & \wh \beta_k\\ && \wh \beta_k & \wh \alpha_k\eem, \qquad \wh{\mbf H}_{k+1,k}=\bem \wh{\mbf H}_k \\ \wh\beta_{k+1}\mbf e_k^\top\eem.$$ The matrix $\wt{\mbf R}_k$ in RSMAR-I is $$\wt{\mbf R}_k=\bem \wh\beta_1\mbf e_1 &\wh{\mbf H}_{k,k-1}\eem=\bem  \wh\beta_1 & \wh \alpha_1 & \wh \beta_2 & & \\ & \wh \beta_2 & \wh\alpha_2 & \ddots &  \\ & &\wh\beta_3 & \ddots & \wh\beta_{k-1}\\  & & & \ddots & \wh \alpha_{k-1}\\ &&&  & \wh \beta_k\eem.$$ 
We need the QR factorization $$\wh{\mbf H}_{k+1,k}=\wh{\mbf Q}_{k+1}\bem\wh{\mbf R}_k\\ \mbf 0\eem, \qquad \wh{\mbf R}_k=\bem \delta_1 & \lambda_1 & \eta_1 &&\\ & \delta_2 & \lambda_2 & \ddots &\\ &&\delta_3& \ddots& \eta_{k-2} \\ && &\ddots & \lambda_{k-1}\\ &&&&\delta_k \eem,$$ where $\wh{\mbf Q}_{k+1}^\top=\mbf G_{k,k+1}\mbf G_{k-1,k}\cdots\mbf G_{1,2}$ is a product of reflections. For $i=1,2,\ldots,k$, the structure of $\mbf G_{i,i+1}$ is $$\mbf G_{i,i+1}=\bem \mbf I_{i-1}  &&&\\ &c_i&s_i& \\ &s_i&-c_i& \\ &&&\mbf I_{k-i}\eem.$$ 
We initialize $\wt\delta_1:=\wh\alpha_1$ and $\wt\lambda_1:=\wh\beta_2$. The $k$th reflection $\mbf G_{k,k+1}$ zeroing out $\wh\beta_{k+1}$ satisfies $$\bem c_k & s_k \\ s_k & -c_k\eem\bem \wt\delta_k & \wt\lambda_k & 0\\ \wh\beta_{k+1}& \wh\alpha_{k+1} & \wh\beta_{k+2} \eem=\bem \delta_k & \lambda_k & \eta_k \\ 0 & \wt\delta_{k+1} & \wt\lambda_{k+1}\eem,$$ where elements  decorated by a tilde are to be updated by the next reflection. Straightforward computations give $$\delta_k:=\sqrt{\wt\delta_k^2+\wh\beta_{k+1}^2},\qquad c_k:=\wt\delta_k/\delta_k, \qquad s_k:=\wh\beta_{k+1}/\delta_k.$$ Then we have the following recursion  \begin{align*}
 \lambda_k := c_k\wt\lambda_k+s_k\wh\alpha_{k+1},\quad 
 \wt\delta_{k+1} := s_k\wt\lambda_k-c_k\wh\alpha_{k+1},\quad 
 \eta_k :=s_k\wh\beta_{k+2},\quad 
 \wt\lambda_{k+1} :=-c_k\wh\beta_{k+2}. 
 \end{align*} 
 The vector $\wh{\mbf t}_{k+1}=\wh{\mbf Q}_{k+1}^\top\wh\beta_1\mbf e_1=\bem \wh t_1 & \wh t_2 & \cdots & \wh t_k & \wt t_{k+1}\eem^\top$ can be obtained by using the recursion $$\wt t_1:=\wh\beta_1, \quad \wh t_i=c_i\wt t_i, \quad \wt t_{i+1}=s_i\wt t_i,\quad i=1,2,\ldots,k.$$ We have $\|\mbf A\mbf r_k\|=|\wt t_{k+1}|=|s_1s_2\cdots s_k\wh\beta_1|.$ By \eqref{imple1}, we have $$\mbf x_k^\rma=\mbf x_0+\bem \mbf r_0 &\wh{\mbf V}_{k-1}\eem\wt{\mbf R}_k^{-1}\wh{\mbf z}_k=\mbf x_0+\bem \mbf r_0 &\wh{\mbf V}_{k-1}\eem\wt{\mbf R}_k^{-1}\wh{\mbf R}_k^{-1}\bem\mbf I_k & \mbf 0\eem\wh{\mbf t}_{k+1}.$$ 
To avoid storing $\wh{\mbf V}_k$, we define 
$$\mbf W_k:=\bem \mbf r_0 &\wh{\mbf V}_{k-1}\eem\wt{\mbf R}_k^{-1}=\bem\mbf w_1 & \mbf w_2 &\cdots & \mbf w_k\eem,\qquad \mbf P_k:=\mbf W_k\wh{\mbf R}_k^{-1}=\bem\mbf p_1 & \mbf p_2 &\cdots & \mbf p_k\eem.$$ Then $$\mbf x_k^\rma=\mbf x_0+\mbf W_k\wh{\mbf R}_k^{-1}\bem\mbf I_k & \mbf 0\eem\wh{\mbf t}_{k+1}=\mbf x_0+\mbf P_k\bem\mbf I_k & \mbf 0\eem\wh{\mbf t}_{k+1}.$$ The columns of $\mbf W_k$ and $\mbf P_k$ can be obtained from the recursions \begin{align*}
&\mbf w_1=\mbf r_0/\wh\beta_1, \quad \mbf w_2=(\wh{\mbf v}_1-\wh\alpha_1\mbf w_1)/\wh\beta_2,\quad \mbf w_k=(\wh{\mbf v}_{k-1}-\wh\beta_{k-1}\mbf w_{k-2}-\wh\alpha_{k-1}\mbf w_{k-1})/\wh\beta_k,\quad k\geq 3,\\ 
&\mbf p_1=\mbf w_1/\delta_1, \quad \mbf p_2=(\mbf w_2-\lambda_1\mbf p_1)/\delta_2,\quad \mbf p_k=(\mbf w_k-\eta_{k-2}\mbf p_{k-2}-\lambda_{k-1}\mbf p_{k-1})/\delta_k,\quad k\geq 3, 	
 \end{align*} and the solution $\mbf x_k^\rma$ may be updated via $$\mbf x_k^\rma=\mbf x_{k-1}^\rma+\wh t_k\mbf p_k.$$
The approach given above is summarized as Algorithm 3. 

\begin{table}[H]
\centering
\begin{tabular*}{170mm}{l}
\toprule {\bf Algorithm 3}. MINARES-I: implementation based on $\mbf A\wh{\mbf V}_k=\wh{\mbf V}_{k+1}\wh{\mbf H}_{k+1,k}$ for $\mcalk_k(\mbf A,\mbf A\mbf r_0)$ \hspace{10mm}
\\ \hline\noalign{\smallskip} \quad {\bf Require}: symmetric $\mbf A\in\mbbr^{n\times n}$, $\mbf b\in\mbbr^n$, $\mbf x_0\in\mbbr^n$, ${\tt tol}>0$, ${\tt maxit}>0$
\\ \noalign{\smallskip} \quad\hspace{.64mm} 1:\ \ $\mbf r_0:=\mbf b-\mbf A\mbf x_0$, $\wh\beta_1:=\|\mbf A\mbf r_0\|$, $\rho_0:=\wh\beta_1$. If $\rho_0<\tt tol$, accept $\mbf x_0^\rma=\mbf x_0$ and exit. 
\\ \noalign{\smallskip} \quad\hspace{.64mm} 2:\ \  $\wh{\mbf v}_1:=\mbf A\mbf r_0/\wh\beta_1$, $\mbf w_1:=\mbf r_0/\wh\beta_1$, $\wt t_1:=\wh\beta_1$, $\wh{\mbf v}_0=\mbf p_{-1}=\mbf p_0=\mbf w_0=\mbf 0$, $c_0=-1$, $s_0=\wt\lambda_0=\eta_{-1}=0$
\\ \noalign{\smallskip} \quad\hspace{.64mm} 3:\ \  {\bf for} $k=1, 2, \ldots, \tt maxit$ {\bf do}
\\ \noalign{\smallskip} \quad\hspace{.64mm} 4:\ \
 \qquad $\wh{\mbf v}_{k+1}:=\mbf A\wh{\mbf v}_k-\wh\beta_k\wh{\mbf v}_{k-1}$ 
\\ \noalign{\smallskip} \quad\hspace{.64mm} 5:\ \ \qquad $\wh\alpha_k:=\wh{\mbf v}_k^\top\wh{\mbf v}_{k+1}$  
\\ \noalign{\smallskip}\hspace{.64mm} \quad 6:\ \ \qquad $\wh\beta_{k+1}\wh{\mbf v}_{k+1}:=\wh{\mbf v}_{k+1}-\wh\alpha_k\wh{\mbf v}_k$   \hfill  {\color[gray]{0.5} $\wh\beta_{k+1}>0$ so that $\|\wh{\mbf v}_{k+1}\|=1$}
\\ \noalign{\smallskip} \quad\hspace{.64mm} 7:\ \  \qquad $\lambda_{k-1}:= c_{k-1}\wt\lambda_{k-1}+s_{k-1}\wh\alpha_k$
\\ \noalign{\smallskip} \quad\hspace{.64mm} 8:\ \ \qquad $\wt \delta_k:=s_{k-1}\wt\lambda_{k-1}-c_{k-1}\wh\alpha_k$
\\ \noalign{\smallskip} \quad\hspace{.64mm} 9:\ \  \qquad $ 
 \eta_{k-1} :=s_{k-1}\wh\beta_{k+1}$
 \\ \noalign{\smallskip} \quad 10:\ \ \qquad $\wt\lambda_k:=-c_{k-1}\wh\beta_{k+1}$
\\ \noalign{\smallskip} \quad 11:\ \ \qquad $\delta_k:=\sqrt{\wt\delta_k^2+\wh\beta_{k+1}^2}$
\\ \noalign{\smallskip} \quad 12:\ \ \qquad $c_k:=\wt\delta_k/\delta_k$
\\ \noalign{\smallskip} \quad 13:\ \  \qquad $s_k:=\wh\beta_{k+1}/\delta_k$
\\ \noalign{\smallskip} \quad 14:\ \  \qquad  $\wh t_k := c_k\wt t_k$
\\ \noalign{\smallskip} \quad 15:\ \  \qquad $\wt t_{k+1}:=s_k\wt t_k$
\\ \noalign{\smallskip} \quad 16:\ \  \qquad $\rho_k:=|\wt t_{k+1}|$
\\ \noalign{\smallskip} \quad 17:\ \  \qquad $\mbf p_k:=(\mbf w_k-\eta_{k-2}\mbf p_{k-2}-\lambda_{k-1}\mbf p_{k-1})/\delta_k$
\\ \noalign{\smallskip} \quad 18:\ \  \qquad $\mbf x_k^\rma:=\mbf x_{k-1}^\rma+\wh t_k\mbf p_k$ 
\\ \noalign{\smallskip} \quad 19:\ \  \qquad {\bf if} $\rho_k<\tt tol$ {\bf then}
\\ \noalign{\smallskip} \quad 20:\ \  \qquad\qquad Accept $\mbf x_k^\rma$ and exit.
\\ \noalign{\smallskip} \quad 21:\ \  \qquad {\bf end if}
\\ \noalign{\smallskip} \quad 22:\ \  \qquad $\mbf w_{k+1}:=(\wh{\mbf v}_k-\wh\beta_k\mbf w_{k-1}-\wh\alpha_k\mbf w_k)/\wh\beta_{k+1}$
\\ \noalign{\smallskip} \quad 23:\ \ {\bf end for} 
\\ 
\bottomrule
\end{tabular*}
\end{table}

\section{Numerical experiments}
We will compare the performance of GMRES, RRGMRES, RSMAR, and DGMRES on singular range-symmetric linear systems, and compare the performance of MINRES-QLP, MINARES, and RSMAR on singular symmetric linear systems. All algorithms stop 
if $k>\tt maxit$, where {\tt maxit} is the maximum number of iterations. In all algorithms, the initial approximate solution $\mbf x_0$ is set to be zero vector. To get a fair comparison, residuals for consistent systems or $\mbf A$-residuals for inconsistent systems are calculated explicitly at each iteration. All experiments are performed using MATLAB R2023b on MacBook Pro with Apple M3 Max chip, 128 GB memory, and macOS Sonoma 14.2.1.

\subsection{Singular range-symmetric linear systems}

In this subsection, we compare the performance of GMRES, RRGMRES, RSMAR, and DGMRES on singular range-symmetric linear systems generated from a matrix arising in the finite difference discretization of the following boundary value problem \begin{align}\label{bvpmodel}\l\{\begin{array}{lcl} \dsp\Delta u+d\frac{\partial u}{\partial x}=f, & \mbox{ in } &\Omega:=[0,1]\times[0,1], \vspace{1mm}\\  u(x,0)=u(x,1), & \mbox{ for }  & 0\leq x\leq 1, \vspace{1mm}\\ u(0,y) =u(1,y), & \mbox{ for }  & 0\leq y\leq 1, \end{array}\r.\end{align} where $d$ is a constant and $f$ is a given function. The matrix $\mbf A$ is given as follows: 
$$\mbf A = \bem \mathbf{T}_m & \mathbf{I}_m & & \mathbf{I}_m \cr \mathbf{I}_m & \ddots & \ddots \cr & \ddots & \ddots & \mathbf{I}_m \cr \mathbf{I}_m & & \mathbf{I}_m & \mathbf{T}_m \eem \in \mathbb{R}^{m^2\times m^2},\quad \mathbf{T}_m = \bem -4 & \alpha_{+} & & \alpha_{-} \cr \alpha_{-} & \ddots & \ddots \cr & \ddots & \ddots & \alpha_{+} \cr \alpha_{+} & & \alpha_{-} & -4 \eem \in \mathbb{R}^{m\times m},$$ where $m=100$, $h=1/m$, $\alpha_\pm=1\pm dh/2$, and $d=10$. This matrix is normal and hence range-symmetric. It is already used to illustrate the performance of GMRES in \cite[Experiment 4.2]{brown1997gmres}. Note that $\mbf A$ is singular with $\nul(\mbf A)=\spa\{\bem 1&1&\cdots &1\eem^\top\}$. 

\begin{figure}[H]
  \centerline{\includegraphics[height=7cm]{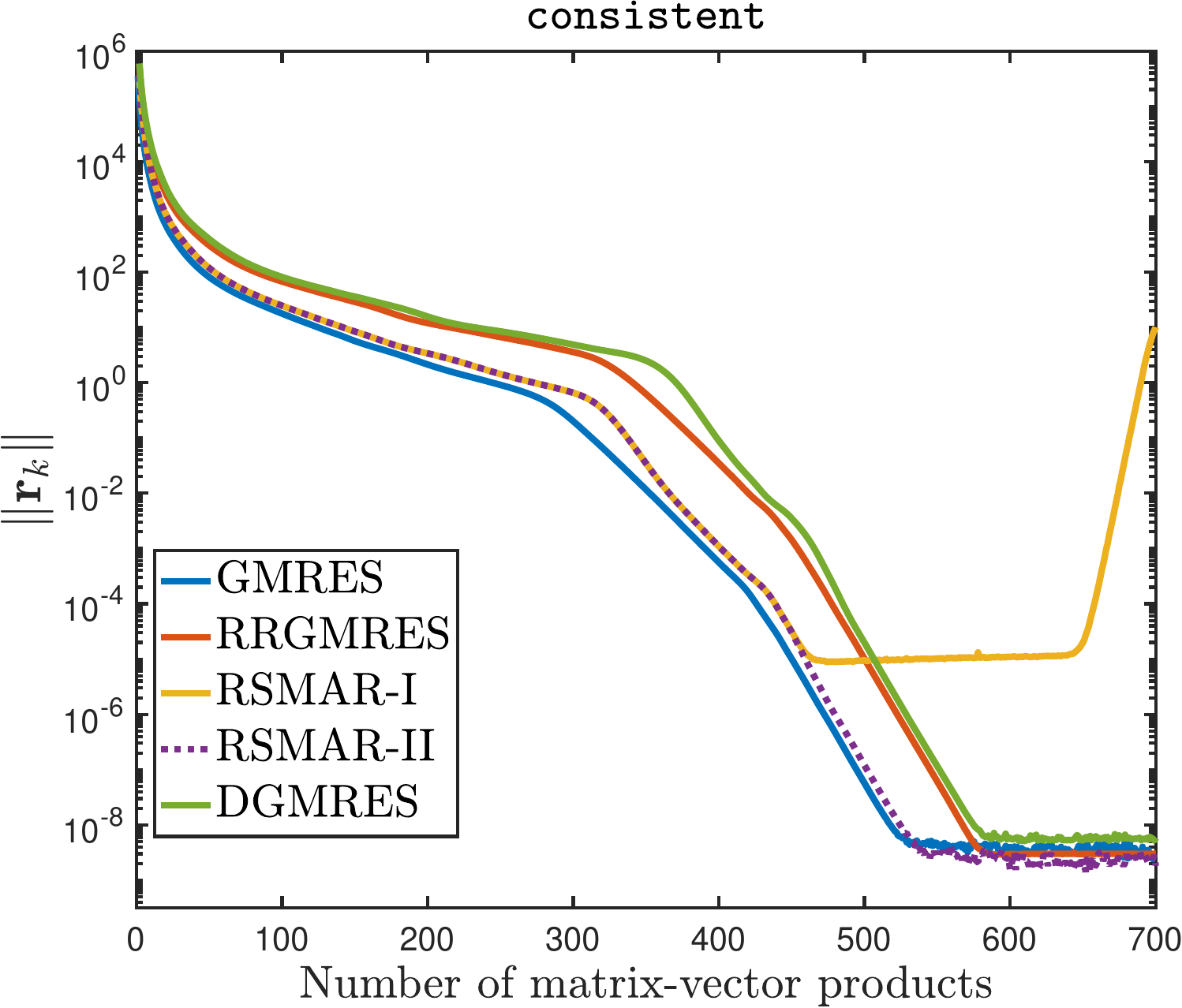}\quad\includegraphics[height=7cm]{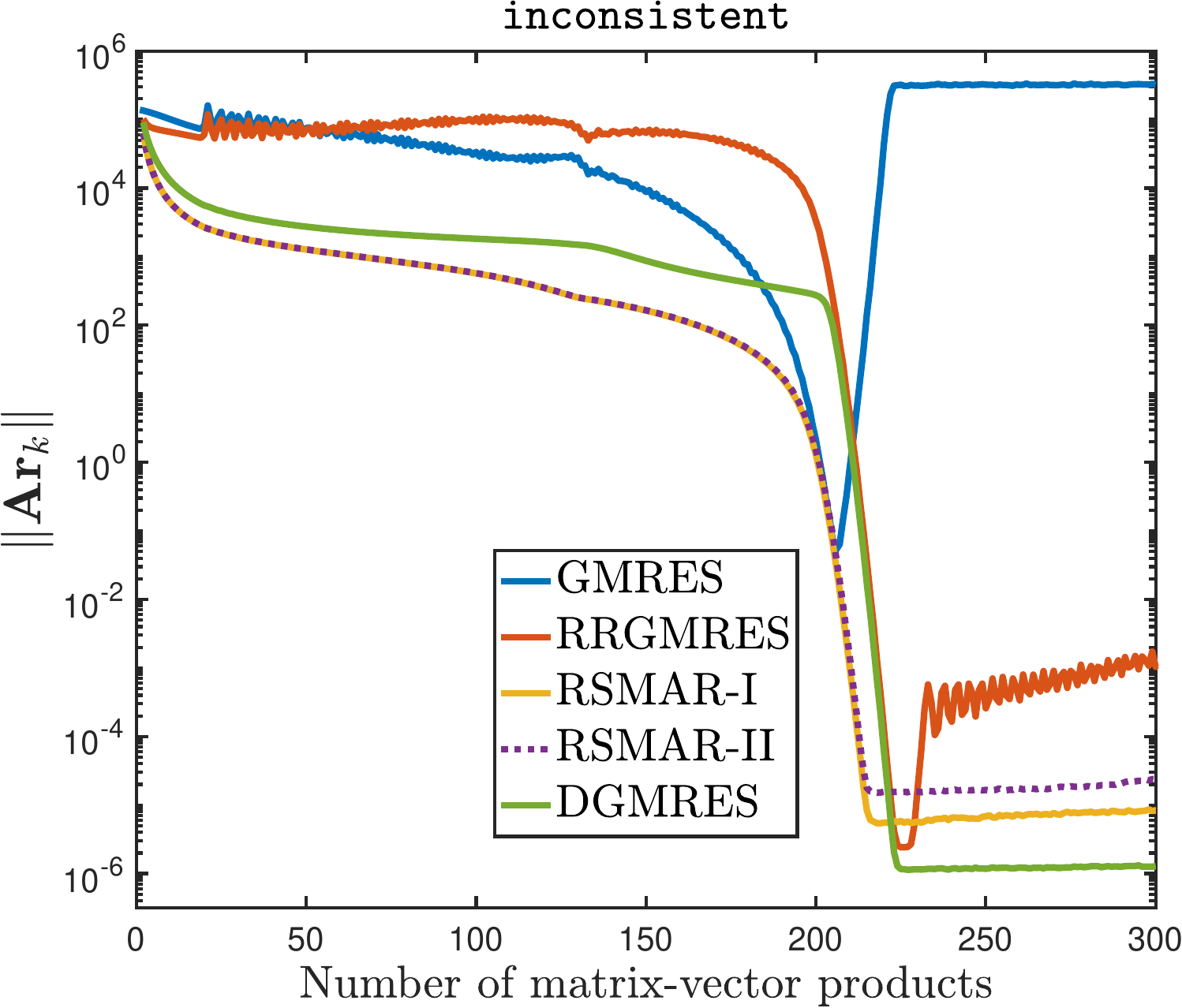}}
  \caption{Residual and $\mbf A$-residual histories for GMRES, RRGMRES, RSMAR, and DGMRES on singular linear systems generated from the matrix arising in the finite difference discretization of the boundary value problem \eqref{bvpmodel}. Left: consistent system with {\tt b = A*rand(m*m,1)}. Right: inconsistent system with $\mbf b$ being a discretization of $f(x,y)=x+y$.}\label{bvp}  
\end{figure}

We first construct a consistent linear system by using MATLAB's script ``{\tt rng("default"); b = A*rand(m*m,1);}''. We then construct an inconsistent linear system by taking $\mbf b$ to be a discretization of $f(x,y)=x+y$. In Figure \ref{bvp}, we plot residual histories for GMRES, RRGMRES, RSMAR, and DGMRES on the consistent system, and $\mbf A$-residual histories for these algorithms on the inconsistent system. We have the following observations. 
\bit \item[(i)] In the consistent case, GMRES, RRGMRES, RSMAR-II, and DGMRES attain almost the same accuracy. RSMAR-I suffers from an instability. The residual norm $\|\mbf r_k\|$ of all algorithms is smooth before reaching the attainable optimal accuracy, and GMRES is slightly faster than other algorithms in terms of number of matrix-vector products.
\item[(ii)] In the inconsistent case, the $\mbf A$-residual norm $\|\mbf A\mbf r_k\|$ of RSMAR and DGMRES is smooth before reaching the attainable optimal accuracy, whereas that of GMRES and RRGMRES is erratic. RSMAR is faster than other algorithms in terms of number of matrix-vector products. The attainable accuracy of DGMRES is the best, and that of GMRES is the worst.
\eit

\subsection{Singular symmetric linear systems}

In this subsection, we compare the performance of MINRES-QLP, MINARES, and RSMAR applied to singular symmetric linear systems generated from symmetric matrices $\mbf A$ from the SuitSparse Matrix Collection \cite{davis2011unive}. Three matrices ({\tt bcsstm36}, {\tt zenios}, and {\tt laser}) are used. In each problem, we scale $\mbf A$ to be $\mbf A/\rho$ with $\rho=\max_{ij}|\mbf A_{ij}|$, so that $\|\mbf A\|\approx 1$. Consistent systems are constructed by using $\mbf b=\mbf A\mbf e$ (with $\mbf e$ a vector of ones), and inconsistent ones are done by using $\mbf b=\mbf e$.  

\begin{figure}[H]
  \centerline{\includegraphics[height=7cm]{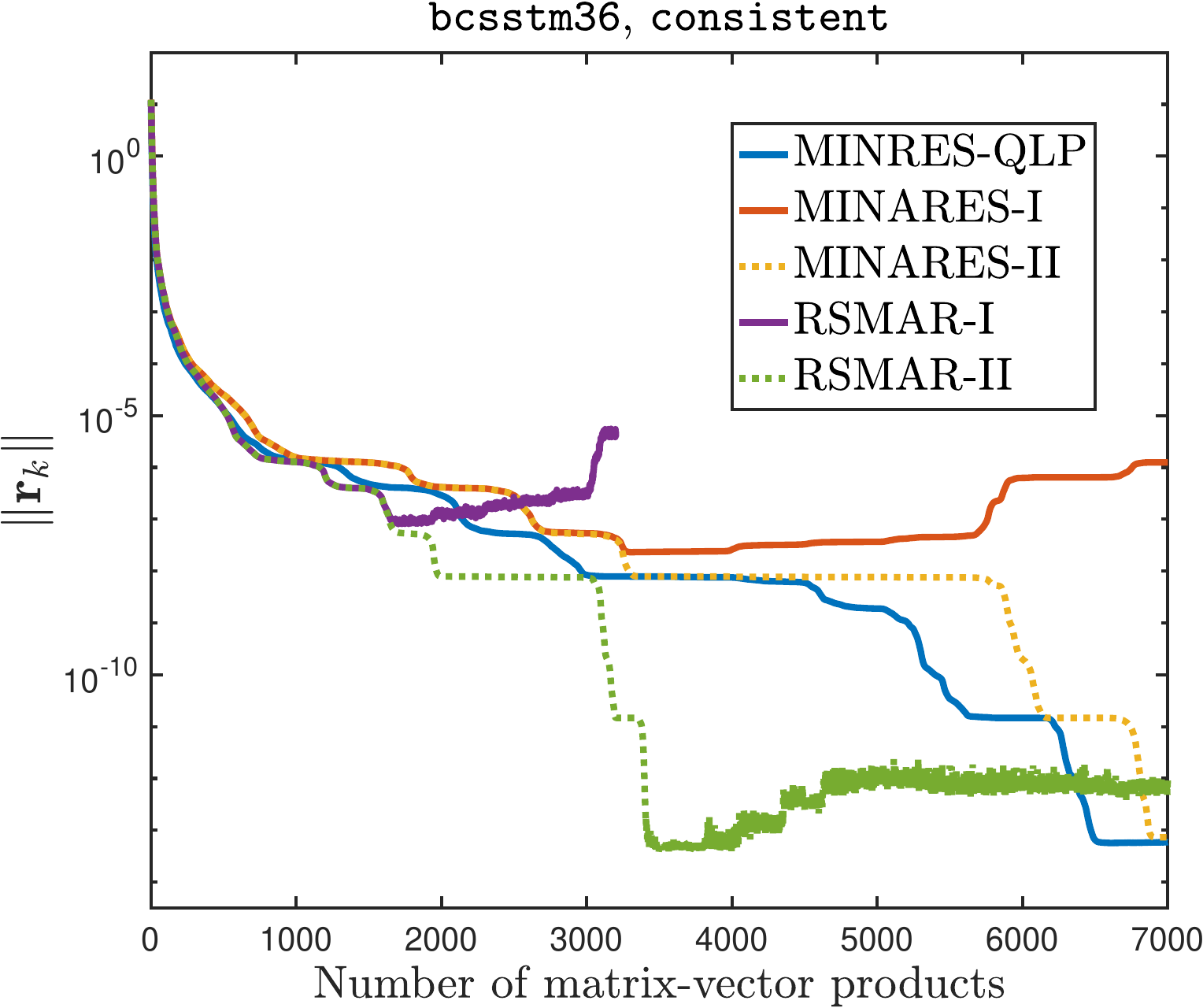}\quad\includegraphics[height=7cm]{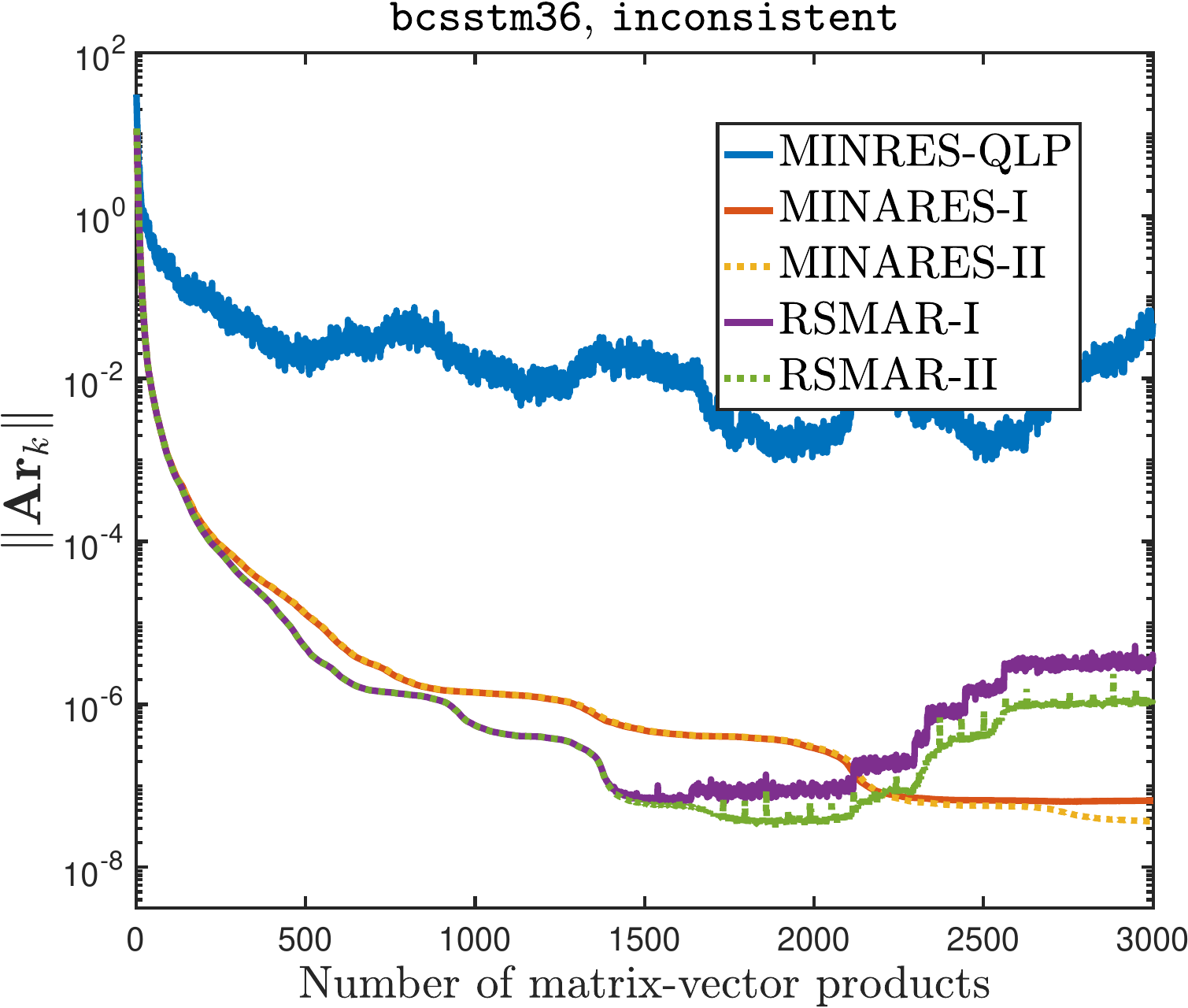}}
  \caption{Residual and $\mbf A$-residual histories for MINRES-QLP, MINARES, and RSMAR on singular linear systems generated from the matrix {\tt bcsstm36} ($n=23052$). Left: consistent system with $\mbf b=\mbf A\mbf e$. Right: inconsistent system with $\mbf b=\mbf e$.}\label{bcsstm36}  
\end{figure}

\begin{figure}[H]
  \centerline{\includegraphics[height=7cm]{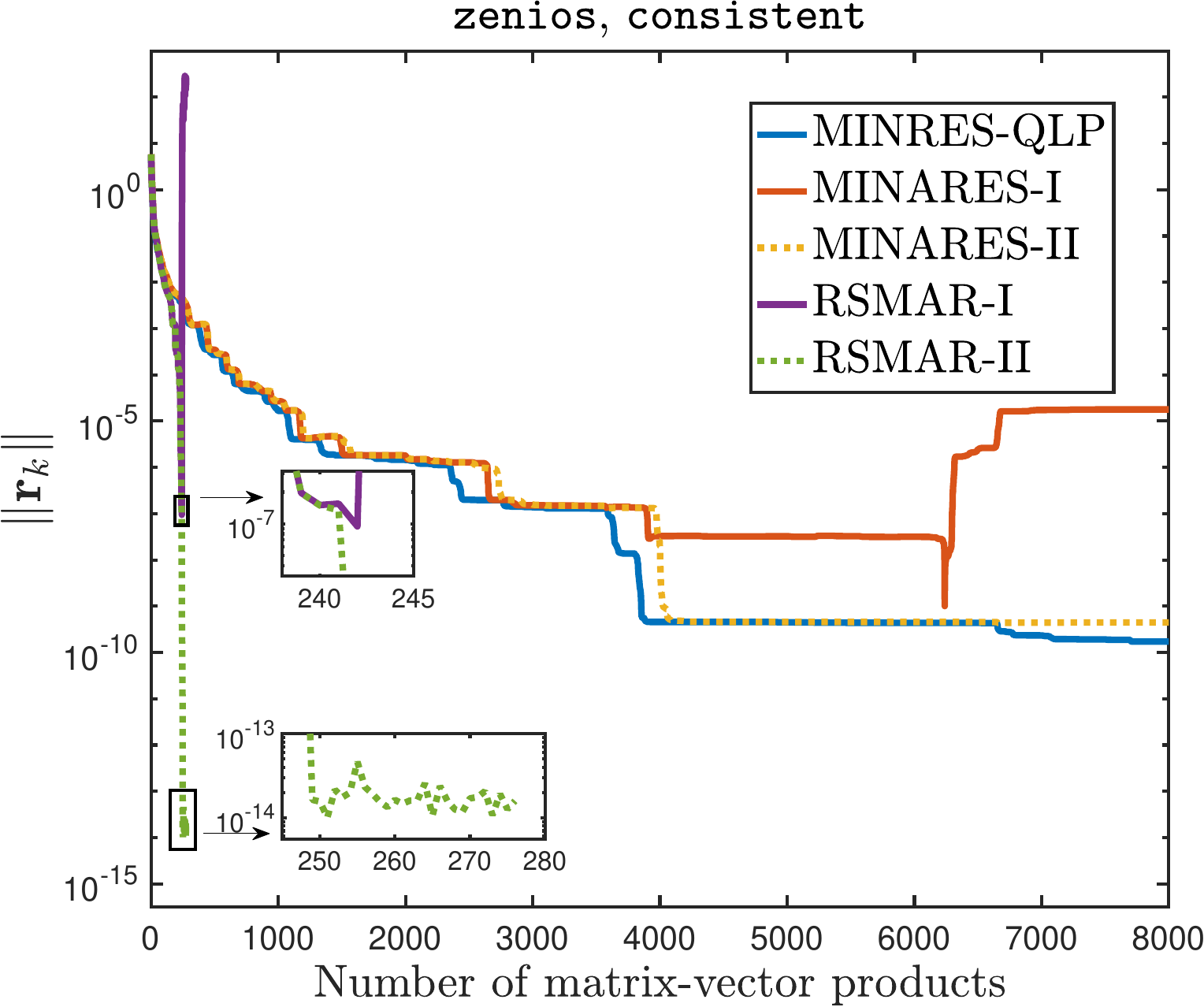}\quad\includegraphics[height=7cm]{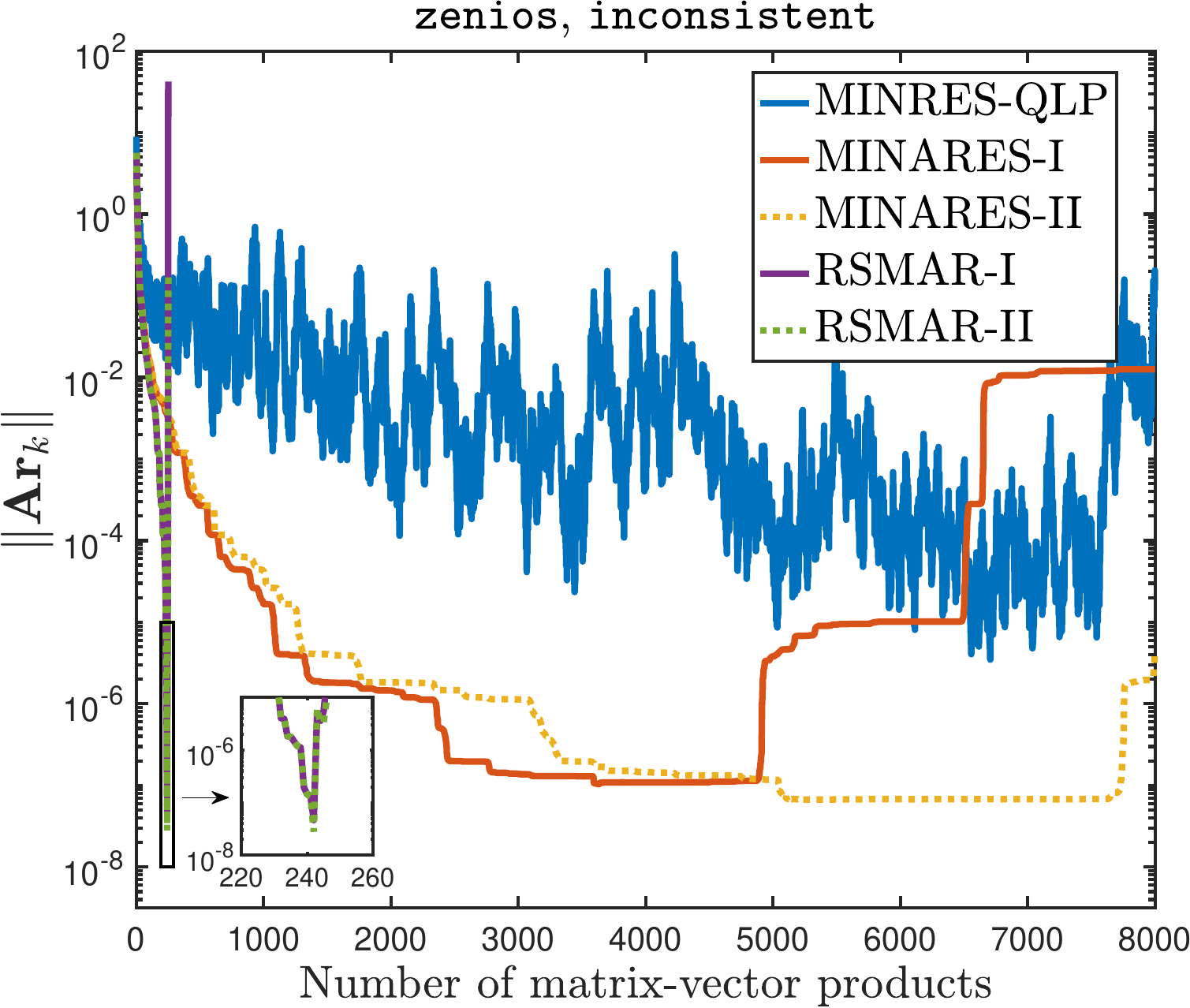}}
  \caption{Residual and $\mbf A$-residual histories for MINRES-QLP, MINARES, and RSMAR on singular linear systems generated from the matrix {\tt zenios} ($n=2873$). Left: consistent system with $\mbf b=\mbf A\mbf e$. Right: inconsistent system with $\mbf b=\mbf e$.}\label{zenios}  
\end{figure}

\begin{figure}[H]
  \centerline{\includegraphics[height=7cm]{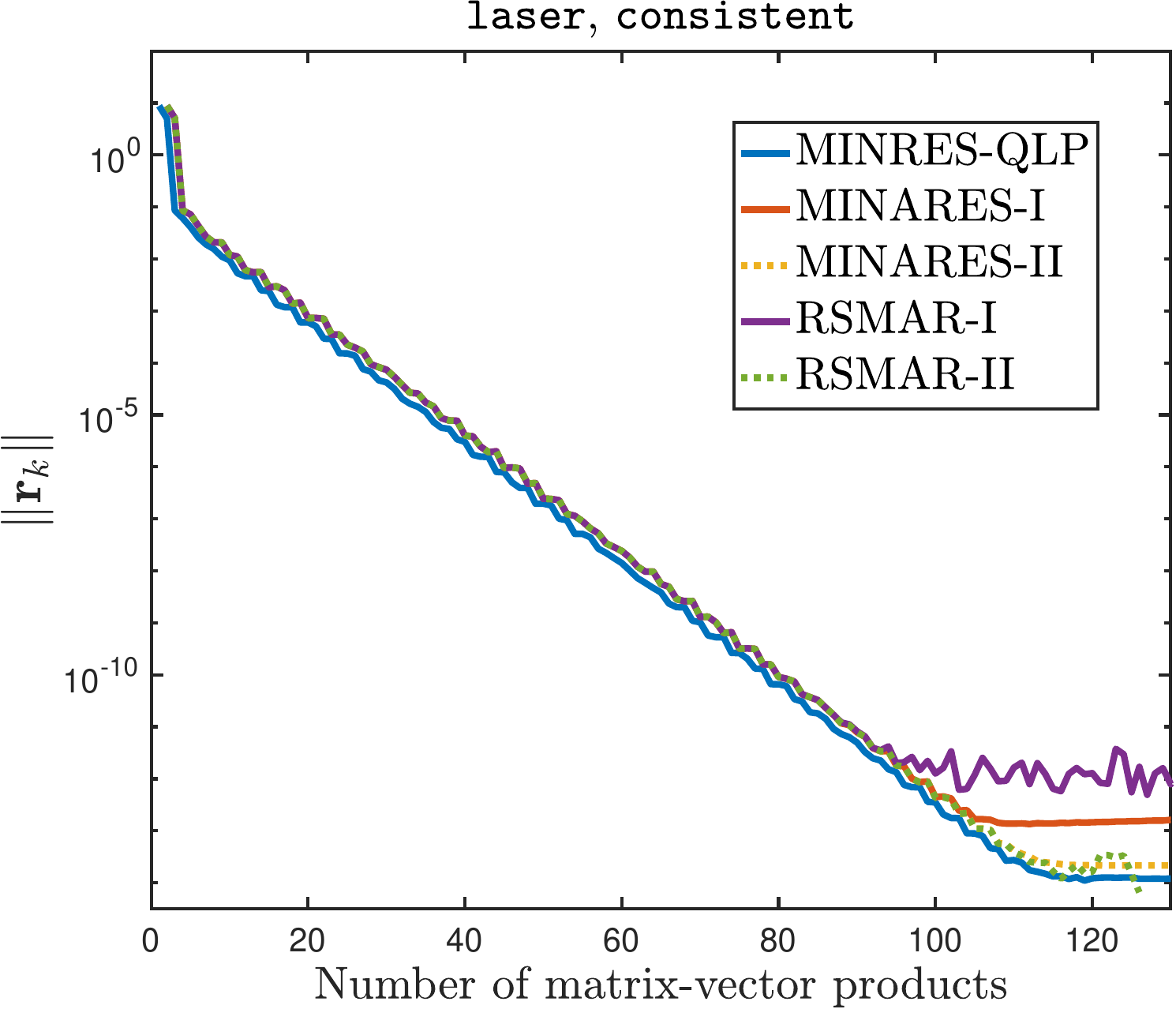}\quad\includegraphics[height=7cm]{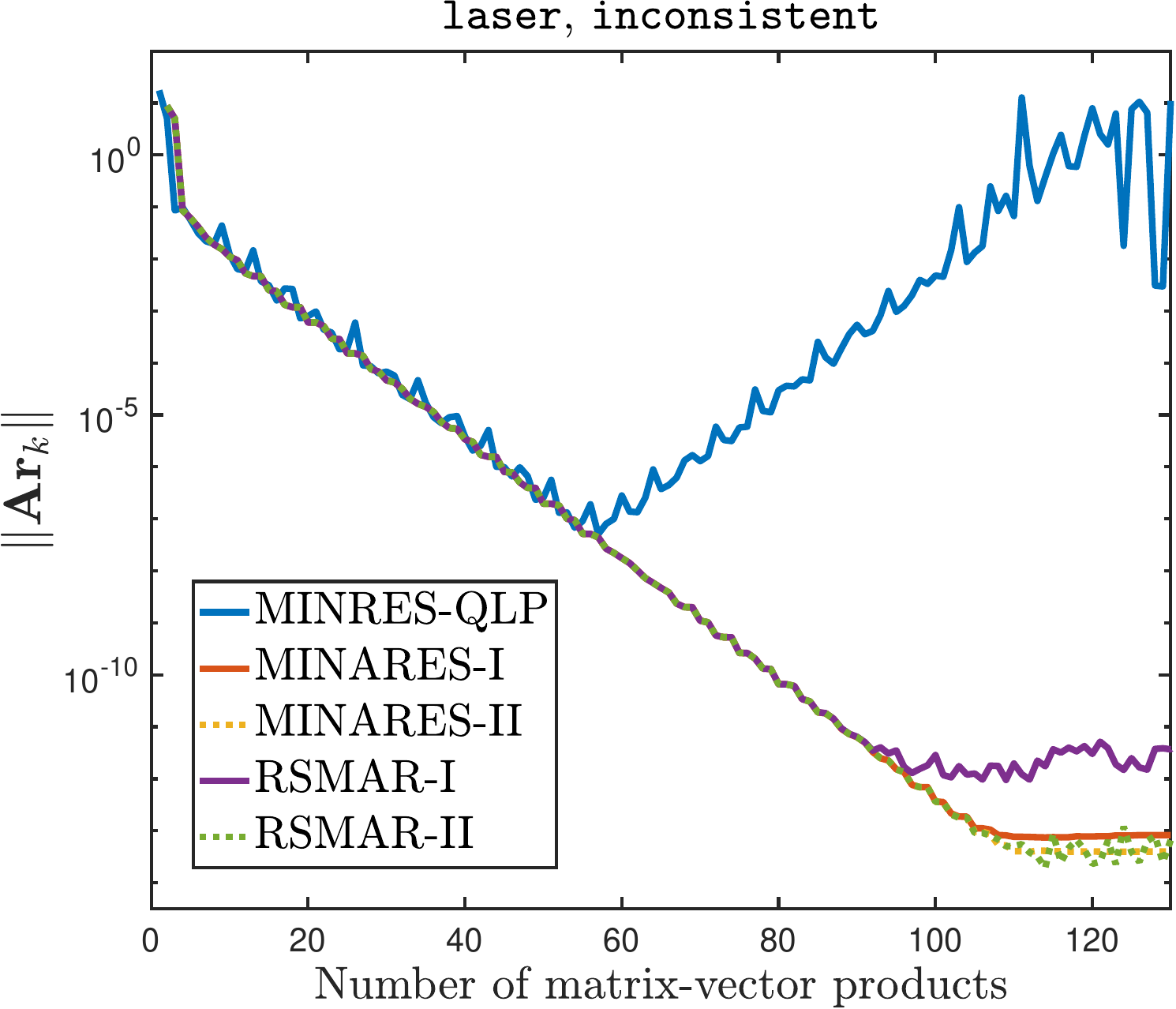}}
  \caption{Residual and $\mbf A$-residual histories for MINRES-QLP and MINARES on singular linear systems generated from the matrix {\tt laser} ($n=3002$). Left: consistent system with $\mbf b=\mbf A\mbf e$. Right: inconsistent system with $\mbf b=\mbf e$.}\label{laser}  
\end{figure}

 We report residual histories for MINRES-QLP, MINARES, and RSMAR on consistent systems, and $\mbf A$-residual histories for these algorithms on inconsistent systems. The MINARES implementation of Montoison, Orban, and Saunders \cite{montoison2023minar} is referred as MINARES-II. Figures \ref{bcsstm36}, \ref{zenios}, and \ref{laser} are on the systems generated using {\tt bcsstm36}, {\tt zenios}, and {\tt laser}, respectively. In the consistent case for the problem {\tt bcsstm36}, RSMAR-I suffers from an instability and we terminate it when the number of iterations $k=3200$. For the problem {\tt zenios}, we terminate RSMAR when $k=275$ in the consistent case, and when $k=250$ in the inconsistent case. We have the following observations. 
 \bit 
 \item[(i)] For all problems, RSMAR-II is better than RSMAR-I and MINARES-II is better than MINARES-I in terms of the attainable optimal accuracy.
 \item[(ii)] For the problems {\tt bcsstm36} and {\tt laser}, RSMAR and MINARES nearly coincide only in the initial phase, and RSMAR-II is faster than MINARES in terms of number of matrix-vector products. For the problem {\tt laser}, RSMAR and MINARES nearly coincide.
 \item[(iii)]  In the consistent cases for the problems {\tt bcsstm36} and {\tt laser}, RSMAR-I suffers from an instability. In all consistent cases, the residual norm $\|\mbf r_k\|$ of all algorithms is smooth before reaching the attainable optimal accuracy. 
 \item[(iv)] In all inconsistent cases, the $\mbf A$-residual norm $\|\mbf A\mbf r_k\|$ of RSMAR and MINARES is smooth before reaching the attainable optimal accuracy, whereas that of MINRES-QLP is erratic. MINRES-QLP suffers from an instability in the inconsistent case for the problem {\tt laser}.
 \eit

\section{Concluding remarks and future work} RSMAR completes the family of Krylov subspace methods based on the Arnoldi process for range-symmetric linear systems. By minimizing the $\mbf A$-residual norm $\|\mbf A\mbf r_k\|$ (which always converges to zero for range-symmetric $\mbf A$), RSMAR can be applied to solve any range-symmetric systems. We have shown that in exact arithmetic, RSMAR and GMRES both determine the pseudoinverse solution if $\mbf b\in\ran(\mbf A)$, and   terminate with the same least squares solution if $\mbf b\notin\ran(\mbf A)$. When the reached least squares solution is not the pseudoinverse solution, the lifting strategy \eqref{lifted} can be used to obtain it. Our numerical experiments show that on singular inconsistent range-symmetric systems, RSMAR outperforms GMRES, RRGMRES, and DGMRES, and should be the preferred method in finite precision arithmetic. As for the implementation for RSMAR, RSMAR-II is better than RSMAR-I in finite precision arithmetic. 

When $\mbf A$ is symmetric, RSMAR is theoretically equivalent to MINARES. The work per iteration and the storage requirements of RSMAR increase with the iterations, while MINARES remains under control even when many iterations are needed. 


There are at least three possible research directions for future work. The first is about preconditioning techniques for RSMAR. The second is about stopping criteria. It would clearly be desirable to terminate the RSMAR iterations when approximately optimal accuracy has been reached. The third is the performance of RSMAR applied to linear discrete ill-posed problems. All of them are being investigated.

Our MATLAB implementations of GMRES, RRGMRES, RSMAR, DGMRES, MINRES-QLP, and MINARES are available at \href{https://kuidu.github.io/code.html}{https://kuidu.github.io/code.html}. The implementations of GMRES, RRGMRES, RSMAR, and DGMRES support restarts. All figures in section 5 can be reproduced by the MATLAB live script {\tt mar.mlx}, which can be obtained from the above website.

%
%
%


	
{\small 

}
\end{document}